\documentclass{siamart190516}

\usepackage[utf8]{inputenc}
\usepackage{url}
\usepackage{xr-hyper}
\usepackage{amsmath,amsfonts, multirow, graphicx}
\usepackage[caption=false]{subfig} 
\usepackage{algorithm,algpseudocode}

\usepackage{xcolor}
 \definecolor{mygrn}{RGB}{0,140,70}

\usepackage[mathscr]{eucal}
\usepackage{amsbsy,bm}
\usepackage{amssymb}
\usepackage{hyperref,cleveref}
\usepackage{tikz}
\usepackage{calc}
\usetikzlibrary{calc,arrows}
\usetikzlibrary{decorations.pathreplacing}
\usepackage{changes}
\usepackage[T1]{fontenc}

\newcommand{\Tra}{\top} 
\newcommand{\Pinv}{\dagger}

\newcommand{\T}[2][]{\boldsymbol{#1\mathscr{\MakeUppercase{#2}}}} 
\newcommand{\M}[2][]{{\bm{#1{\mathbf{\MakeUppercase{#2}}}}}}

\newcommand{\Tentry}[1]{\ensuremath\T{#1}(i_1,\dots,i_N)}
\newcommand{\Mz}[3][]{\M[#1]{#2}_{(#3)}} 

\newcommand{\TT}[2]{\T{T}_{\T{#1},#2}} 
\newcommand{\slice}[2]{} 
\newcommand{\TTslice}[2]{ 
\ifx1#2  
	\renewcommand{\slice}[2]{\TT{#1}{1}(i_1,:)}
\else \ifx#2N 
	\renewcommand{\slice}[2]{\TT{#1}{N}(:,i_N)}
\else 
	\renewcommand{\slice}[2]{\TT{#1}{#2}(:,i_{#2},:)}
\fi \fi
\slice{#1}{#2} 
} 
\newcommand{\HOp}{\mathcal{H}}
\newcommand{\VOp}{\mathcal{V}}

\newcommand{\R}{\mathbb{R}}
\newcommand{\mc}[1]{\mathcal{#1}}

\makeatletter
\newcommand*{\addFileDependency}[1]{
  \typeout{(#1)}
  \@addtofilelist{#1}
  \IfFileExists{#1}{}{\typeout{No file #1.}}
}
\makeatother

\title{Randomized algorithms for rounding in the Tensor-Train format}

\author{Hussam Al Daas\thanks{Computational Mathematics Group, Rutherford Appleton Laboratory, UK. (\email{hussam.al-daas@stfc.ac.uk})} \and Grey Ballard\thanks{Department of Computer Science, Wake Forest University, USA. (\email{ballard@wfu.edu})} \and Paul Cazeaux\thanks{Department of Mathematics, University of Kansas, USA. (\email{pcazeaux@ku.edu}, \email{amiedlar@ku.edu})} \and Eric Hallman\thanks{Department of Mathematics, North Carolina State University, USA. (\email{erhallma@ncsu.edu}, \email{twreid@ncsu.edu},  \email{asaibab@ncsu.edu})}\and Agnieszka Mi\k{e}dlar\footnotemark[3]  \and Mirjeta Pasha\thanks{School of Mathematical and Statistical Sciences, Arizona State University, USA. (\email{mpasha3@asu.edu})} \and Tim W.\ Reid\footnotemark[4] \and Arvind K.\ Saibaba\footnotemark[4]}

\begin{document}
\maketitle

\begin{abstract}
The Tensor-Train (TT) format is a highly compact low-rank representation for high-dimensional tensors. TT is particularly useful when representing approximations to the solutions of certain types of parametrized partial differential equations. For many of these problems, computing the solution explicitly would require an infeasible amount of memory and computational time. While the TT format makes these problems tractable, iterative techniques for solving the PDEs must be adapted to perform arithmetic while maintaining the implicit structure. The fundamental operation used to maintain feasible memory and computational time is called \emph{rounding}, which truncates the internal ranks of a tensor already in TT format. We propose several randomized algorithms for this task that are generalizations of randomized low-rank matrix approximation algorithms and provide significant reduction in computation compared to deterministic TT-rounding algorithms. Randomization is particularly effective in the case of rounding a sum of TT-tensors (where we observe $20\times$ speedup), which is the bottleneck computation in the adaptation of GMRES to vectors in TT format. We present the randomized algorithms and compare their empirical accuracy and computational time with deterministic alternatives.
\end{abstract}


\section{Introduction}
An increasing number of applications in science and technology involve the manipulation of multi-dimensional data, or tensors that are higher order equivalents of vectors (first-order) and matrices (second-order). The number of elements of a tensor as well as the storage consumption grow exponentially with the number of the dimensions, a phenomenon  known as the \emph{curse of dimensionality}.  When problems of high dimensions are concerned, beating the {curse of dimensionality} and finding a solution efficiently remains a challenge. Nevertheless, different tensor formats and methods based on tensor products \cite{khoromskij2012tensors, kolda2009tensor, grasedyck2013literature, hackbusch2012tensor, oseledets2012solution} have shown potential for mitigating the curse of dimensionality and tackling high-dimensional problems that could not be addressed with conventional methods. Initially, the concept of tensor decompositions was introduced in 1927 by expressing a tensor as the sum of a finite number of rank-one tensors \cite{hitchcock1927expression} --- also known as the \emph{canonical format}. The canonical format's memory requirements are not high, though it can suffer from numerical stability issues \cite{de2008tensor, holtz2012alternating}.
Tensors in Tucker form \cite{beck2000multiconfiguration} are well known in quantum chemistry \cite{de2008tensor, holtz2012alternating} since they yield robust algorithms due to the ability to form an embedded manifold \cite{koch2010dynamical}, but one of the disadvantages of the Tucker format is its storage consumption that still depends exponentially on the number of dimensions.

One of the most promising tensor formats is the Tensor-Train (TT) format, a tensor product format that was initially proposed in quantum physics, also known as \emph{matrix product states} (MPS) \cite{fannes1992finitely}, and was reinvented in numerical linear algebra \cite{oseledets2009breaking, oseledets2011tensor}.
It combines both the advantages of the canonical and Tucker formats, i.e., 1) the storage consumption of a tensor depends linearly on the number of dimensions and 2) there
exist robust algorithms for the computation of best approximations.
Applications of the TT format arise from various applications such as high-dimensional PDEs like the Fokker Planck equations \cite{dolgov2012fast, richter2021solving}, quantum physics \cite{schollwock2005density, meyer2009basic}, high-dimensional data analysis \cite{klus2015numerical, klus2018tensor}, machine learning \cite{beylkin2009multivariate, duvenaud2015advances, cohen2016expressive, obukhov2021spectral}, and uncertainty quantification \cite{zhang2014enabling, loukrezis2018high} to mention just a few. Typically, those applications require an approximate solution of linear systems of equations, eigenvalue problems, or completion problems \cite{gelss2017nearest, ballani2013projection, rauhut2015tensor}. The TT format is a low-rank representation that, for TT-tensors with small rank, offers a tremendous reduction in the computational complexity and often exposes the structure of the problem. The use of low-rank structures such as the 
TT format \cite{oseledets2011tensor} to represent
high-dimensional objects allows the solution of linear high-dimensional problems by generalizing standard numerical linear algebra techniques to multi-index arrays of coefficients (tensors) and the multivariate functions they approximate. 

In this paper, we focus on the problem of rounding a tensor in TT format; that is, assuming that we are given a TT-tensor, we want to find a compressed representation that is nearly as accurate as the original representation. There are several techniques for computing the initial TT-tensors which do not require forming the entire tensor explicitly~\cite{DolS20,KreKNT15,OseT10,SavO11}. One such technique that is popularly used is called the TT-cross approximation. The standard TT approach to rounding, proposed by Oseledets~\cite{oseledets2011tensor},  has two phases~\cite[Algorithm 2]{oseledets2011tensor}: orthogonalization followed by compression (typically using the SVD). Here, by orthogonalization, we mean a sweep of orthogonalization steps across every tensor core. Analysis shows that the orthogonalization step dominates the computational cost of this approach. Motivated by this observation, the goal of this work is to develop randomized algorithms for rounding TT-tensors that avoid expensive orthogonalization.
In the following, we present the main contributions of this paper.

\paragraph{Overview of the paper and main contributions} 

This paper develops several new randomized algorithms for rounding tensors in the TT format and is organized as follows. In \Cref{sec: background}, we set some notation as well as review some basic material on randomized matrix algorithms and standard TT operations along with a detailed analysis of their computational costs. In \Cref{sec:randttalgs}, we propose various new randomized algorithms for TT-rounding with the focus on Randomize-then-Orthogonalize, Two-Sided-Randomization, and rounding of a sum of TT-tensors.
\begin{enumerate}
    \item In \cref{alg:TT-rounding-Rand1}, Orthogonalize-then-Randomize, we replace the SVD step in the standard TT-rounding algorithm with a randomized SVD assuming that the truncated ranks are known \textit{a priori}.       \item In \Cref{alg:TT-rounding-Rand2}, Randomize-then-Orthogonalize, we propose to form randomized sketches of each core by nested contractions with a TT-tensor with random cores in a first step, before performing the orthogonalization sweep on these much smaller matrices. Our analysis and experiments show that this approach allowed for the best speedup compared to the deterministic algorithm while retaining excellent accuracy.
    \item In \Cref{alg:TT-rounding-Rand3}, Two-Sided-Randomization, we completely eliminate the need for separate orthogonalization and compression sweeps. Instead, we work with a two-sided randomized approach which computes products with two random tensors followed by a compression step (which involves orthogonalization of much smaller matrices). {Although this approach is slightly more expensive in terms of flops count and less accurate than techniques mentioned before, it eliminates the need of extensive orthogonalization and allows for the truncation phase to be more independent and highly parallelizable.}
    \item We extend the Randomize-then-Orthogonalize approach for compressing a TT-tensor that is presented as a sum of TT-tensors (\Cref{alg:TT-rounding-sum-Rand2}). This special case is of importance in many applications such as solving parametric linear systems in the TT format. The use of randomization enables significant performance improvements by exploiting the structure of the sum tensor in a way that a deterministic algorithm cannot. 
\end{enumerate}
We provide an analysis of the computational cost of the proposed algorithms in \Cref{ssec:compcosts}  and show that they are computationally more efficient than existing algorithms. We justify our analysis through numerical experiments in \Cref{sec:numex} on both synthetic data and tensors generated while solving parametric partial differential equations (PDEs). Some conclusions and future outlook are presented in \Cref{sec: conclusions}. The {\sc Matlab} code for the implementation and numerical experiments is publicly available at \href{https://github.com/SAMSI-RandTensors/randomizedTT}{https://github.com/SAMSI-RandTensors/randomizedTT}.

\paragraph{Related work} 
There have been several recent developments in obtaining low-rank compression of tensors. We limit our literature review to the publications dealing with TT-tensors described in \cite{oseledets2011tensor}, which is closest to our work, and refer the reader to review papers for other developments in tensor decompositions~\cite{ahmadi2021randomized,cichocki2016tensor,cichocki2017tensor,grasedyck2013literature}. Oseledets~\cite{oseledets2011tensor} proposes a method for rounding TT-tensors. A parallel version of this method is introduced and developed in \cite{daas2020parallel}. Our newly proposed approaches are more computationally efficient compared to existing deterministic algorithms. Other works~\cite{huber2017randomized,che2019randomized,ahmadi2020randomized} discuss randomized algorithms for compressing tensors in the TT format. These approaches differ from ours in that they require access to the entries of the tensor, i.e., they do not assume that the tensor is already in TT format. 
A recent paper~\cite{alger2020tensor} also uses randomization to produce a TT approximation of a full tensor but relies on tensor actions (i.e., applications of the tensor on $N-1$ vectors, where $N$ is the order of the tensor). 
Other methods for constructing a low-rank compression in the TT format involve alternating least squares~\cite{dolgov2014alternating}. 
The use of tensor random projections in which the random tensors are taken to be in TT format have also been considered in~\cite{batselier2018computing,rakhshan2020tensorized,feng2020tensor}. 
While these papers use randomization in the context of TT-tensors, none of them directly address the problem of rounding which is the central focus of our paper.

\section{Background}
\label{sec: background}
Here, we review the notation and necessary operations involving tensors in a modest amount of detail. For a more comprehensive exposition we refer the reader to~\cite{oseledets2011tensor,kolda2009tensor,daas2020parallel}.

\subsection{Notation}

We denote tensors by boldface script letters (e.g., $\T{X}$) and matrices by boldface Roman letters (e.g., $\M{A}$). 
We follow {\sc Matlab}-like convention and denote the entries of a 3-way tensor $\T{X}$ as $\T{X}(i,j,k)$. 
A colon denotes the entire range of indices in that dimension. 
We denote the column fibers as $\T{X}(:,j,k)$, row fibers as $\T{X}(j,:,k)$ and tube fibers as $\T{X}(j,k,:)$.
The \emph{mode-$n$ unfolding} (or matricization) of the tensor $\T{X} $ is denoted as $\M{X}_{(n)} \in \mathbb{R}^{I_n \times ({I}/{I_n})}$, where $I = I_1I_2\cdots I_n$. The columns of the mode-$n$ unfolding are composed of the appropriate mode-$n$ fibers, e.g., the columns of mode-$1$ unfolding are column fibers and the columns of mode-$3$ unfolding are tube fibers. Given a matrix $\M{A} \in \mathbb{R}^{ M\times I_n}$, the mode-$n$ product $\T{Y} = \T{X}\times_n\M{A}$ is defined by its mode-$n$ unfolding $\M{Y}_{(n)} = \M{A}\M{X}_{(n)}$.
The norm of a tensor is equivalent to the Frobenius norm of any of its unfoldings: $\|\T{X}\|=\|\M{X}_{(n)}\|_F$.

An \emph{order-$N$ tensor} $\T{X} \in \mathbb{R}^{I_1 \times \cdots \times I_N}$ is in the TT format if there exist positive integers $R_0, \ldots, R_N$ with $R_0 = R_N = 1$ and order-3 tensors $\TT{X}{1}, \ldots, \TT{X}{N} $, called \emph{TT-cores}, with $\TT{X}{n} \in \mathbb{R}^{R_{n-1} \times I_n \times R_n}$ for $1 \leq n \leq N$, such that
\[
	\Tentry{X} = \TTslice{X}{1} \cdot \ldots \cdot \TTslice{X}{n} \cdot \ldots \cdot \TTslice{X}{N},
\]
where $1\leq i_n\leq I_n$.
Note that because $R_0=R_N=1$, the first and last TT-cores are (order-2) matrices so $\TTslice{X}{1} \in \mathbb{R}^{R_1}$ and $\TTslice{X}{N} \in \mathbb{R}^{R_{N-1}}$.
The $R_{n-1} \times R_n$ matrix $\TTslice{X}{n}$ is referred to as the $i_n$th slice of the $n$th TT-core of $\T{X}$. It is worth mentioning that the TT decomposition is not unique due to the multiplicative nature of the format. 
\begin{figure}[h!]
 \input{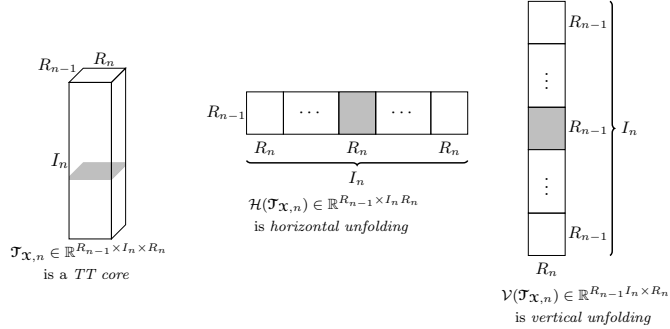}
    \centering
    \scalebox{0.7}{

  \begin{tabular}{ccc}
    \begin{minipage}{0.23\linewidth}
      \begin{tikzpicture}[scale=.7,textnode/.style={scale=.5}]
        \pgfmathsetmacro{\ione}{.25*\none}
        \pgfmathsetmacro{\itwo}{.6*\ntwo}
        \fill[gray!50] (0,-\itwo,0) rectangle (-\rfive,-\itwo-\vecwidth,0);
        \shade[top color=gray!50,bottom color=gray!50] (0,-\itwo,0) -- ++(0,0,-\rtwo) -- ++(0,-\vecwidth,0) -- ++(0,0,\rtwo) -- cycle;
        \shade[top color=gray!50,bottom color=gray!50] (0,-\itwo,0) -- ++(0,0,-\rtwo) -- ++(0,-\vecwidth,0) -- ++(0,0,\rtwo) -- cycle;
        \shade[top color=gray!50,bottom color=gray!50] (0,-\itwo,0) -- ++(-\rfive,0,0) -- ++(0,0,-\rtwo) -- ++(\rfive,0,0) -- cycle;
        \drawTower
      \end{tikzpicture}\\
      \centering{\footnotesize 
        $\TT{X}{n} \in \R^{\nthcoredimone \times \nthcoredimtwo \times \nthcoredimthree}$ \\ 
        is a \emph{TT core} 
      }
    \end{minipage}&
    \begin{minipage}{0.42\linewidth}
      \begin{tikzpicture}[scale=.7,textnode/.style={scale=.5}]
        \drawHUnfolding
      \end{tikzpicture}
      \centering{\footnotesize 
        $\HOp(\TT{X}{n}) \in \R^{\nthcoredimone \times \nthcoredimtwo\nthcoredimthree}$ \\  is \emph{horizontal unfolding}
      }
    \end{minipage}&
    \begin{minipage}{0.25\linewidth}
      \begin{tikzpicture}[scale=.7,textnode/.style={scale=.5}]
        \centering
        \drawVUnfolding
      \end{tikzpicture}
      \centering{\footnotesize 
        $\VOp(\TT{X}{n}) \in \R^{\nthcoredimone\nthcoredimtwo \times \nthcoredimthree}$ \\ is \emph{vertical unfolding}
      }
    \end{minipage}
  \end{tabular}
    }
    \caption{Horizontal and vertical unfoldings of a TT-core $\TT{X}{n}$.}
    \label{fig:tt_unfoldings}
\end{figure}

\begin{figure}[!ht]
    \centering
    \scalebox{0.6}{
\tikzstyle{vertex} = [circle,fill=blue!40,draw,scale=.8]
\tikzstyle{edge} = [line width=1pt]

\begin{tikzpicture}

\pgfmathsetmacro{\d}{5} 
\pgfmathsetmacro{\dm}{int(\d-1)}

\foreach \i in {1,...,\d} {
    \node[vertex] (a\i) at (\i*2,0) {$\TT{X}{\i}$};
    \node (b\i) at (\i*2,-1.5) {};
}
\foreach \i in {1,...,\dm} {
    \pgfmathsetmacro{\j}{int(\i+1)}
    \draw[edge] (a\i) -- (a\j) node[midway,above] {$R_{\i}$};
}
\foreach \i in {1,...,\d} {
    \draw[edge] (a\i) -- (b\i) node[midway,right] {$I_{\i}$};
}

\end{tikzpicture}
    }
    \caption{Tensor network diagram for an order-$5$ TT-tensor.}
    \label{fig:tt_network}
\end{figure}
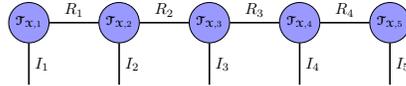

In order to express the arithmetic operations on TT-cores using linear algebra, we will often use two specific matrix unfoldings of the order-$3$ tensors.
The \emph{horizontal unfolding} of a TT-core $\TT{X}{n}$ corresponds to the concatenation of the slices $\TTslice{X}{n}$ for $i_n = 1,\ldots, I_n$ horizontally.
We denote the corresponding operator by $\HOp$, so that $\HOp(\TT{X}{n})$ is an $R_{n-1} \times I_n R_n$ matrix.
The \emph{vertical unfolding} of a TT-core $\TT{X}{n}$ corresponds to the concatenation of the slices $\TTslice{X}{n}$ for $i_n = 1,\ldots, I_n$ vertically.
We denote the corresponding operator by $\VOp$, so that $\VOp(\TT{X}{n})$ is an $R_{n-1} I_n \times R_n$ matrix,  see \Cref{fig:tt_unfoldings}. Moreover, we will often make use of  a \textit{tensor network diagram}, see \Cref{fig:tt_network}, to graphically illustrate TT-tensor operations.
Here nodes represent tensors and edges represent modes so that connected nodes can be contracted.

Let $\M{X}_{(1:n)} \in \mathbb{R}^{(I_1I_2\cdots I_n)\times (I_{n+1}\cdots I_N)}$ denote an unfolding of the first $n$ modes of a TT-tensor $\T{X}$. It has the rank $R_n$ representation
\[\M{X}_{(1:n)} = \VOp(\TT{X}{1:n})\HOp(\TT{X}{n+1:N}),\]
where in an extension of their earlier definitions, $\VOp(\TT{X}{1:n})\in \mathbb{R}^{(I_1I_2\cdots I_n)\times R_n}$ represents the mode-$(n+1)$ unfolding of the product of the first $n$ TT-cores and $\HOp(\TT{X}{n+1:N})\in \mathbb{R}^{R_n \times (I_{n+1}\cdots I_N)})$ represents the mode-1 unfolding of the product of the final $N-n$ TT-cores.
Likewise, we can write the same unfolding as a product of four matrices, see~\cite[Eq.~(2.3)]{daas2020parallel}, i.e.,
\begin{equation}
    \label{eq:four_prod_mat_rep}
\M{X}_{(1:n)} = (\M{I}_{I_n} \otimes \VOp(\TT{X}{1:n-1}))\VOp(\TT{X}{n})\HOp(\TT{X}{n+1})(\HOp(\TT{X}{n+2:N}) \otimes \M{I}_{I_{n+1}}).
\end{equation}

Suppose we have two tensors $\T{Y}$ and $\T{Z}$ of the same dimension, and consider their sum $\T{X}$. The cores of the tensor $\T{X}$ can be expressed as 
\[ \TT{X}{n}(:,i_n,:) = \begin{bmatrix} \TT{Y}{n}(:,i_n,:) & \\ & \TT{Z}{n}(:,i_n,:) \end{bmatrix} \qquad 2\leq n \leq N-1, \]
and for the first and the last core, we have
\[ \TT{X}{1}(i_1, :) = \begin{bmatrix} \TT{Y}{1}(i_1,:) & \TT{Z}{1}(i_1,:) \end{bmatrix} \quad \mbox{ and } \quad \TT{X}{N}(:,i_N) =\begin{bmatrix} \TT{Y}{N}(:,i_N)\\ \TT{Z}{N}(:,i_N)  \end{bmatrix}. \]

Let $\M{X} \in \mathbb{R}^{m\times n}$ with $m \geq n$. We denote the \emph{thin QR} factorization of $\M{X}$ as $\M{X} = \M{QR}$, where $\M{Q} \in \R^{m \times n}$ has orthonormal columns and $\M{R} \in \R^{n \times n}$ is upper triangular; we also write $[\M{Q},\M{R}] = \text{QR}(\M{X})$ for use in algorithms. The SVD of $\M{X}$ is denoted by $\M{U\Sigma V}^\top$, where the matrix $\M{U}\in \mathbb{R}^{m\times n}$ has orthonormal columns containing the left singular vectors, $\M\Sigma \in \mathbb{R}^{n\times n}$ is a diagonal matrix with the singular values on the diagonal and $\M{V} \in \mathbb{R}^{n\times n}$ is an orthogonal matrix, whose columns contain the right singular vectors. Assuming that the Householder QR algorithm is used and $\M{Q}$ is formed explicitly, the computational cost of the QR factorization is $4mn^2 - \frac{4n^3}{3} +  \mc{O}(n^2)$ flops.  Given a threshold $\varepsilon > 0$, we truncate the singular values of $\M{X}$ to obtain a rank-$k$ approximation $\M{U}_k\M\Sigma_k\M{V}_k^\top$ of matrix $\M{X}$, which satisfies $\|\M{X}-\M{U}_k\M\Sigma_k\M{V}_k^\top \|_F \leq \varepsilon \|\M{X}\|_F $. This is denoted as $[\M{U}_k,\M\Sigma_k,\M{V}_k] =$\; SVD$(\M{X},\varepsilon)$. The computational cost of computing the SVD is $\mc{O}(mn^2)$ flops.

\subsection{Randomized matrix algorithms}\label{ssec:randsvd} 
An important component of our approach is the use of randomized matrix methods for low-rank matrix approximation. In this subsection, we briefly review a few well-established randomized algorithms. 

The first algorithm is the basic version of the randomized SVD proposed in~\cite{halko2011finding}. Suppose we want to compute a low-rank approximation of a matrix $\M{X} \in \mathbb{R}^{m\times n}$; let the target rank be denoted by $r$ and we pick an oversampling parameter $p$ such that $r + p \leq \min\{m,n\}$. We generate a random matrix $\M{\Omega} \in \mathbb{R}^{n\times (r+p)}$; in practice, we take the entries of this matrix to be i.i.d.~standard Gaussian random variables. Then, we compute the product $\M{Y} = \M{X\Omega}$ and obtain its thin QR factorization $\M{Y} = \M{QR}$. The main insight exploited by randomized SVD is that if the rank of $\M{X}$ is close to $r$, or the singular values of $\M{X}$ decay rapidly beyond $r$, then the range of $\M{Q}$ approximates well the range of $\M{X}$ in the sense that $\M{X} \approx \M{QQ}^\top\M{X}$; we then use $\M{QQ}^\top\M{X}$ as a low-rank approximation to $\M{X}$. The computational cost of this approach is 
\[ C_\text{randSVD} = 2(r+p)mn + \mathcal{O}(r^2(m+n)) \quad \text{flops}. \]
Additional postprocessing can be performed to convert the low-rank approximation in the SVD format, or to truncate the low-rank approximation to rank$-r$;  see~\cite{halko2011finding} for additional details.

There is one variant of this algorithm that is of particular importance to our newly proposed methods, {the generalized Nystr\"{o}m method~\cite{nakatsukasa2020fast}}. 
The generalized Nystr\"{o}m method avoids the orthogonalization step {when computing a low-rank approximation} by using a two-sided randomized approach. Let us define two Gaussian random matrices $\M{\Omega}\in\mathbb{R}^{n\times s}$ and $\M{\Psi}\in\mathbb{R}^{t\times m}$, where $ r \leq s \leq \min\{m,n\}$ (note that $t$ also satisfies a similar inequality). 
A low-rank approximation to $\M{X}$ is computed as 
\begin{equation}
\label{eq:EssenNY}
\M{X} \approx \M{Y} (\M{\Psi}\M{X\Omega})^\dagger \M{Z},
\end{equation}
where $\M{Y} = \M{X\Omega}$ and $\M{Z} = \M{\Psi X}$. To implement the pseudoinverse,~\cite{nakatsukasa2020fast} suggests computing the QR factorization $\M{\Psi}\M{X\Omega} = \M{QR}$ and then obtaining the low-rank approximation $(\M{YR}^{-1}) (\M{Q}^\top \M{Z})$.  
If the low-rank approximation is desired in the SVD format, this can be done by additional post-processing.   In~\cite{nakatsukasa2020fast}, the author recommends setting the sketch parameters as $s = r$ and $t=\lceil 1.5 r \rceil$. The associated computational cost is 
\[
C_\text{genNys} = 2mn(s+t) + \mathcal{O}(t^2(m+n) + ts^2) \quad  \mbox{flops}. 
\]

\subsection{Standard TT arithmetic}\label{ssec:ttrounding}
In this subsection, we review the standard approach to TT-rounding, first proposed in~\cite{oseledets2011tensor}, using the notation of~\cite{daas2020parallel}. We also review the concepts of tensor contractions.

To explain the rounding procedure for TT-tensors, we consider the following analogy from matrices. Let $\M{Y} = \M{AB}$ be an outer product matrix where $\M{A}$ is $m\times r$ and $\M{B}$ is $r\times n$ and $r \leq \min\{m,n\}$. To obtain an approximation of $\M{Y}$ with rank $\ell < r$, we employ an orthogonalization step followed by a compression step. In the orthogonalization step, we want to make $\M{Y}$ right orthogonal. That is, we compute the thin QR factorization $\M{B}^\top = \M{QR}$, and then compute $\M{Z} = \M{A}\M{R}^\top$. This gives 
$\M{Y} = \M{AB} = \M{Z}\M{Q}^\top$,
where $\M{Q}^\top$ has orthonormal rows. In the second step, we compress $\M{Z}$ by computing the rank-$\ell$ truncated SVD $\M{Z} \approx \M{U}\M{\Sigma} \M{V}_Z^\top$. To obtain an overall low-rank approximation to $\M{Y}$, we compute $\M{V} = \M{QV}_Z $, so that $\M{Y} \approx \M{ U \Sigma V}^\top$. 

Following~\cite{daas2020parallel}, we say a tensor is \textit{right orthogonal} if its horizontal unfoldings $\HOp(\TT{X}{n})$ have orthonormal rows for $n=2,\dots,N$ (all except the first core).
Similarly, we say that a tensor is \textit{left orthogonal} if its vertical unfoldings $\VOp(\TT{X}{n})$ have orthonormal columns for $n=1,\dots,N-1$ (all except the last core).

\paragraph{Right-to-Left Orthogonalization} 

Suppose we are given a TT-tensor $\T{Y}$. To obtain a right orthogonal TT-tensor $\T{X}$ equivalent to $\T{Y}$, we first compute the thin QR factorization $\M{Q}\M{R} = \HOp(\TT{Y}{N})^\top$ and set the core tensors $\TT{X}{N-1}$ and $\TT{X}{N}$ as 
\[
\VOp(\TT{Y}{N-1})\HOp(\TT{Y}{N}) = \VOp(\TT{Y}{N-1})(\M{Q}\M{R})^\top = \underbrace{(\VOp(\TT{Y}{N-1})\M{R}^\top)}_{\VOp(\TT{X}{N-1})} \underbrace{(\M{Q}^\top)}_{\HOp(\TT{X}{N})}. 
\]
This procedure is continued through cores $N-1,\dots,2$ but we do not orthogonalize the first core. The details of right-to-left orthogonalization are given in \Cref{alg:orthogonalizationr2l} which will form the foundation for many of the subsequent algorithms. We can similarly obtain a left orthogonal tensor by processing the modes starting from mode-1, but we omit the details here.

\begin{algorithm}[!ht]
 \caption{Right-to-Left Orthogonalization}
\label{alg:orthogonalizationr2l}
\begin{algorithmic}[1]
    \Require{A tensor $\T{Y}$ in TT format}
    \Ensure{$\T{X}$ is right-orthogonal tensor equivalent to $\T{Y}$}
    \Function{$\T{X} = $ OrthogonalizeRL}{$\T{Y}$} 
    \State {$\TT{X}{N}=\TT{Y}{N}$}
	\For{$n=N$ down to $2$}
		\State {$[\HOp(\TT{X}{n})^\top,\M{R}] = \text{QR}(\HOp(\TT{X}{n})^\top)$  \Comment{thin QR factorization} \label{line:orthogonalizationr2l:qr}
		\State $\VOp(\TT{X}{n-1}) = \VOp(\TT{Y}{n-1}) \cdot \M{R}^\top$ \Comment{$\TT{X}{n-1} = \TT{Y}{n-1} \times_3 \M{R}^\top $}} \label{line:orthogonalizationr2l:matmul}
	\EndFor
	\EndFunction
\end{algorithmic}
\end{algorithm}

\paragraph{TT-Rounding} 

Suppose, now, that we want to round the tensor $\T{Y}$ in the TT format, i.e., compress the TT format of a tensor by decreasing the TT-ranks $\{R_n\}$. In the first step of the TT-rounding approach, we first obtain a tensor $\T{X}$ that is right-orthogonal by applying \Cref{alg:orthogonalizationr2l}. Starting with mode-$1$, for each mode, we compute a low-rank approximation of the vertical unfolding $\VOp(\TT{X}{n})$; rather than computing an SVD directly, we first compute the thin QR factorization of $\VOp(\TT{X}{n})$; followed by an SVD of the upper triangular factor $\M{R}$. We then obtain a low-rank approximation 
$\VOp(\TT{X}{n}) \approx \M[\widehat]{U}\M[\widehat]{\Sigma}\M[\widehat]{V}^\top.$
The number of singular values and vectors retained in the low-rank approximation, depend on the threshold $\varepsilon = \frac{\|\T{Y}\|}{\sqrt{N-1}}\varepsilon_0$, where $\varepsilon_0$ is a user-defined threshold that controls the overall accuracy. We then rewrite $\VOp(\TT{X}{n})$ by combining it with the low-rank factor as $\VOp(\TT{X}{n}) = \VOp(\TT{X}{n})\M[\widehat]{U}$. The other two factors $\M[\widehat]{\Sigma}\M[\widehat]{V}^\top$ are combined with the horizontal unfolding $\HOp(\TT{X}{n+1})$ for processing at the next step. This process is terminated after $N-1$ steps and the resulting tensor $\T{X}$ satisfies $\|\T{X}-\T{Y}\| \leq \varepsilon_0 \|\T{Y}\|.$ The details are given in \Cref{alg:TT-rounding-RLR}. 

\begin{algorithm}[!ht]
\caption{TT-Rounding}
\label{alg:TT-rounding-RLR}
\begin{algorithmic}[1]
	\Require{A tensor $\T{Y}$ in TT format, user-defined threshold $\varepsilon_0 > 0$ } 
	\Ensure{A tensor $\T{X}$ in TT format with reduced ranks such that $\|\T{X} - \T{Y}\| \leq \varepsilon_0 \|\T{Y}\|$}
	\Function{$\T{X} =$ TT-Rounding}{$\T{Y}, \varepsilon_0$}
	    \State{ $\T{X} = $ \Call{OrthogonalizeRL}{$\T{Y}$} } \label{line:TTRound:OrthRL}
		\State{Compute $\|\T{Y}\|_F $  and the truncation threshold $\varepsilon = \frac{\|\T{Y}\|}{\sqrt{N-1}} \varepsilon_0$} \label{line:TTRound:Line3}
		\State{ Set $\TT{X}{1} = \TT{Y}{1}$.}
		\For{$n=1$ to $N-1$} \label{line:TTRound:loop}
		    \State{$[\VOp(\TT{X}{n}),\M{R}] = \textsc{QR}(\VOp(\TT{X}{n}))$ \Comment{thin QR factorization}} \label{line:TTRound:loopQR}
			\State{$ [\M[\widehat]{U}, \M[\widehat]{\Sigma}, \M[\widehat]{V}] = \textsc{SVD}(\M{R},\varepsilon)$ \Comment{$\varepsilon$-truncated SVD factorization}} \label{line:TTRound:loopSVD}
			\State{$\VOp(\TT{X}{n}) = \VOp(\TT{X}{n}) \M[\widehat]{U}$ \Comment{$\TT{X}{n}=\TT{X}{n} \times_3 \M[\widehat]{U}$}} \label{line:TTRound:loopMatMul1}
			\State{$\HOp(\TT{X}{n+1}) = \M[\widehat]{\Sigma} \M[\widehat]{V}^\Tra \HOp(\TT{X}{n+1})$ \Comment{$\TT{X}{n+1}=\TT{X}{n+1} \times_1 (\M[\widehat]{\Sigma}\M[\widehat]{V}^\Tra)$}} \label{line:TTRound:loopMatMul2}
		\EndFor \label{line:TTRound:Line10}
	\EndFunction
\end{algorithmic}
\end{algorithm}

\paragraph{Right-to-Left Partial Contraction}

We consider two TT-tensors $\T{X}$ and $\T{Y}$ with ranks $\{ R^{\T{X}}_j\}$ and $\{ R^{\T{Y}}_j\}$, respectively. 
For $n = 2, \ldots, N$ we define the partial contraction matrices
\begin{equation}\label{eqn:partial_contraction}
	\M{W}_{n-1} = \HOp(\TT{X}{n:N})\HOp(\TT{Y}{n:N})^\top \in \R^{R_{n-1}^{\T{X}}\times R_{n-1}^{\T{Y}}}.
\end{equation}
These partial contractions can be computed sequentially as 
\begin{equation}\label{eq:RLContraction}
	\begin{split}
		\VOp(\TT{Z}{n}) &= \VOp(\TT{X}{n})\M{W}_n, \\ 
		\M{W}_{n-1} &= \HOp(\TT{Z}{n})\HOp(\TT{Y}{n})^\top
	\end{split}
\end{equation}
for ${n = 2, \ldots, N-1}$, with $\M{W}_{N-1} = \HOp(\TT{X}{N})\HOp(\TT{Y}{N})^\top$. 
Here $\T{Z}$ is a temporary TT-tensor with compatible dimensions and ranks. 

\begin{figure}
    \input{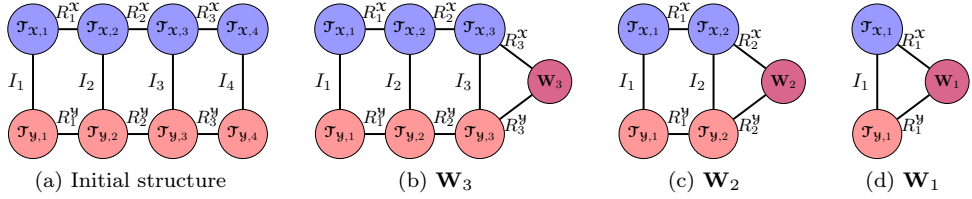}
    \newcommand{\pcscale}{.62}
     \hspace{-0.2in}
    \subfloat[Initial structure]{
        \label{fig:partial_contraction_a}
        \scalebox{0.75}{
        \begin{tikzpicture}[scale=\pcscale]
            \drawcores{4} 
        \end{tikzpicture}
        }
    }
    \centering
    \subfloat[$\M{W}_{3}$]{
        \label{fig:partial_contraction_b}
        \scalebox{0.75}{
        \begin{tikzpicture}[scale=\pcscale]
            \drawcores{3} 
            \drawW{3} 
        \end{tikzpicture}
        }
    }
    \centering
    \subfloat[$\M{W}_2$]{
        \label{fig:partial_contraction_c}
        \scalebox{0.75}{
        \begin{tikzpicture}[scale=\pcscale]
            \drawcores{2} 
            \drawW{2} 
        \end{tikzpicture}
        }
    }
    \centering
    \subfloat[$\M{W}_{1}$]{
        \label{fig:partial_contraction_d}
        \scalebox{0.75}{
        \begin{tikzpicture}[scale=\pcscale]
            \drawcores{1} 
            \drawW{1} 
        \end{tikzpicture}
        }
    }
    \caption{Right-to-Left partial contraction steps for $N=4$.}
    \label{fig:partial_contraction}
\end{figure}

The process of computing the matrices $\{\M{W}_{n-1}\}_{n=2}^N$ according to \cref{eq:RLContraction} is called
a \emph{right-to-left partial contraction} of tensors $\T{X}$ and $\T{Y}$, and is illustrated in \Cref{fig:partial_contraction}. The corresponding algorithm is presented in \Cref{alg:contractionr2l}.

\begin{algorithm}[!ht]
 \caption{Right-to-Left Contraction of Tensors $\T{X}$ and $\T{Y}$. }
 \label{alg:contractionr2l}
\begin{algorithmic}[1]
    \Require{Tensors $\T{X},\T{Y}$ with consistent dimensions in TT format and ranks $\{ R^{\T{X}}_n\}$ and $\{ R^{\T{Y}}_n\}$, respectively. }
    \Ensure{Matrices $\{\M{W}_n\}$ satisfy $\M{W}_n = \HOp(\TT{X}{n{+}1:N})\HOp(\TT{Y}{n{+}1:N})^\top$ for $1 \leq n < N$}
    \Function{$[\{\M{W}_n\}] = $ PartialContractionsRL}{$\T{X},\T{Y}$}
    \State{$\M{W}_{N-1} = \HOp(\TT{X}{N}) \HOp(\TT{Y}{N})^\Tra$} \label{line:contractionr2l:MatMul0}
		\For{$n=N-1$ down to $2$}
			\State{$\VOp(\TT{Z}{n}) =  \VOp(\TT{X}{n}) \M{W}_{n}$ \Comment{$\TT{Z}{n} = \TT{X}{n} \times_3 \M{W}_n$, for temporary $\TT{Z}{n}$}} \label{line:contractionr2l:MatMul1}
			\State{$\M{W}_{n-1} = \HOp(\TT{Z}{n})\HOp(\TT{Y}{n})^\Tra$ \Comment{matrix multiplication, $\M{W}_{n-1}$ is $R^{\T{X}}_{n-1} \times R^{\T{Y}}_{n-1}$}} \label{line:contractionr2l:MatMul2}
		\EndFor
	\EndFunction
\end{algorithmic}
\end{algorithm}

Detailed analysis of the overall computational costs of \emph{Right-to-Left Orthogonalization} (\Cref{alg:orthogonalizationr2l}), \emph{TT-Rounding} (\Cref{alg:TT-rounding-RLR}) and \emph{Right-to-Left Contraction} (\Cref{alg:contractionr2l}) is presented in \Cref{sec:CompCostStandard}.

%
%

\section{Randomized Algorithms for TT Rounding}\label{sec:randttalgs}
In this section, we propose three new randomized algorithms to perform rounding of a tensor in the TT format, i.e., given an original TT-tensor $\T{Y}$ with TT-ranks $\{R_n\}$ we seek a compressed TT-tensor representation $\T{X}$ with \textit{a priori} known target ranks $\{\ell_n\}$. In randomized SVD, see~\Cref{ssec:randsvd}, it is common to include an oversampling term; that is, if we seek a rank-$r$ decomposition of a matrix $\M{X}$, we use the number of samples (alternatively, columns of $\M{\Omega}$) as $\ell = r+p$, where $r$ is the target rank and $p$ is the oversampling parameter. The resulting low-rank approximation $\M{QQ}^\top\M{X}$ is of rank $\ell$. However, in the TT case, to save on notation, when we say target TT-ranks $\{\ell_n\}$, we assume that this rank automatically includes the necessary oversampling parameter.

\subsection{Orthogonalize-then-Randomize}

The first algorithm we propose is very similar to the standard TT-rounding algorithm; the main difference is that we replace the truncated SVD step in \Cref{alg:TT-rounding-RLR} with the basic version of the randomized SVD reviewed in \Cref{ssec:randsvd}. The nomenclature of this algorithm is clear from the fact that there are two phases in this approach: an (already discussed) orthogonalization phase followed by a compression phase which utilizes randomized SVD. 
%
\begin{algorithm}[!ht]
\caption{TT-Rounding: Orthogonalize-then-Randomize}
\label{alg:TT-rounding-Rand1}
\begin{algorithmic}[1]
	\Require{A tensor $\T{Y}$ in TT format with ranks $\{R_n\}$, target TT-ranks $\{\ell_n\}$}
	\Ensure{A tensor $\T{X}$ in TT format with ranks $\{\ell_n\}$}
	\Function{$\T{X} =$ TT-Rounding-OrthRand}{$\T{Y}, \{\ell_n\}$}
	    \State $\T{Y} = $ \Call{OrthogonalizeRL}{$\T{Y}$} \label{line:Rand1:orth}
		\For{$n=1$ to $N-1$} \label{line:Rand1:for}
		    \State{$\M{Z}_n = \VOp(\TT{Y}{n})$  \Comment{$\TT{Y}{n}$ is $\ell_{n-1} \times I_n \times R_n$}} 
		    \State{$\M{Y}_n = \M{Z}_n \M{\Omega}_n$ \Comment{form the sketched matrix}} \label{line:Rand1:sketch}
			\State{$[\VOp(\TT{X}{n}),\sim] = \textsc{QR}(\M{Y}_n)$ \Comment{thin QR to compute an orthonormal basis}} \label{line:Rand1:qr}
			\State{$\M{M}_n = \VOp(\TT{X}{n})^\Tra \M{Z}_n$ \Comment{form $\ell_n \times R_n$ matrix}} \label{line:Rand1:proj1}
			\State{$\HOp(\TT{X}{n+1}) = \M{M}_n \HOp(\TT{Y}{n+1})$ \Comment{ $\TT{X}{n+1} = \TT{Y}{n+1} \times_1 \M{M}_n$}} \label{line:Rand1:proj2}
		\EndFor \label{line:Rand1:endfor}
	\EndFunction
\end{algorithmic}
\end{algorithm}

In the first version, we assume that given a tensor $\T{Y}$ with TT-ranks $\{R_n\}$ we want to obtain a compressed representation with ranks $\{\ell_n\}$; that is, we assume that the target TT-ranks are known in advance. The first phase, i.e., the orthogonalization phase is accomplished using \Cref{alg:orthogonalizationr2l} to obtain the tensor $\T{Y}$ which is right-orthogonal and equivalent to $\T{Y}$. In the second phase, we loop over the cores of the tensor $\T{Y}$ (excluding the last core): for each core, we apply randomized SVD to the vertical unfolding $\M{Z}_n = \VOp(\TT{Y}{n})$. That is, we compute $\M{Q}_n$ such that 
\[  \M{Z}_n \approx \M{Q}_n\M{Q}_n^\top \M{Z}_n, \]
by first generating a random Gaussian matrix $\M\Omega_n \in \mathbb{R}^{R_n\times \ell_n}$. We then compute $\M{Y}_n = \M{Z}_n\M\Omega_n$ and its thin QR factorization to obtain an orthonormal basis $\M{Q}_n$ for the range of $\M{Y}_n$. Finally, we set $\VOp(\TT{X}{n}) = \M{Q}_n$ and $\HOp(\TT{X}{n+1}) = (\M{Q}_n^\top \M{Z}_n) \HOp(\TT{Y}{n+1})$ to obtain the compressed tensor $\T{X}$. The details  are given in \Cref{alg:TT-rounding-Rand1}.

We also investigated a version of \cref{alg:TT-rounding-Rand1} that does not require the ranks of the rounded tensors to be known in advance.  In this approach, we replace \Cref{line:Rand1:sketch,line:Rand1:qr} of \Cref{alg:TT-rounding-Rand1} with a randomized range finder algorithm \cite{martinsson2016randomized,yu2018efficient}. This adaptive method would produce a tensor $\T{X}$ that satisfies the desired tolerance $\|\T{X} - \T{Y}\| \leq \varepsilon_0 \|\T{Y}\|$, where $\varepsilon_0$ is a user defined threshold. While this approach produced TT-tensors with the desired tolerance, in numerical experiments it did not yield significant speedup over deterministic TT-rounding (\cref{alg:TT-rounding-RLR}). Some insight into this behavior is given in \Cref{ssec:compcosts} which shows that the orthogonalization step is the dominant computational cost of deterministic and Orthogonalize-then-Randomize algorithms. Therefore, we did not pursue this approach further.

\subsection{Randomize-then-Orthogonalize}\label{ssec:rto}
To motivate the next algorithm, we consider the overall computational cost of \Cref{alg:TT-rounding-Rand1}, see \Cref{ssec:compcosts}, which is dominated by the first, orthogonalization, phase of the algorithm. First, we consider a new Randomize-then-Orthogonalize algorithm that uses randomization to reduce the overall computation cost of the TT-rounding procedure. It works by avoiding an expensive orthogonalization of the original TT-tensor $\T{Y}$ with TT-ranks $\{R_n\}$ and instead uses randomization to reduce the computational cost. In contrast to the next approach in \Cref{ssec:twosided} (Two-Sided-Randomization), here we use randomization only on one side.

We first offer a way to construct random Gaussian TT-tensors whose cores are composed of independent random Gaussian entries.
\begin{definition}[Random Gaussian TT-Tensor]\label{def:randtt}
Given a set of target TT-ranks $\{\ell_n\}$, we generate a random Gaussian TT-tensor $\T{R} \in \mathbb{R}^{I_1\times \dots I_N}$ such that each  core tensor $\TT{R}{n} \in \mathbb{R}^{\ell_{n-1}\times I_n\times \ell_n}$ is filled with random, independent, normally  distributed entries with mean $0$ and variance $1/(\ell_{n-1} I_n \ell_n)$ for $1\leq n\leq N$.
\end{definition}

By this definition, while the cores of $\T{R}$ have independent entries, the entries of the full tensor themselves are not independent. This normalization is chosen such that $\mathbb{E} \|\TT{R}{n}\|_F^2 = 1$ and is sometimes necessary to ensure that no overflow occurs during the rounding computations. Note that constructing this random tensor requires only generating and storing $\sum_{n=1}^N \ell_{n-1}\ell_nI_n$ random entries. A related but distinct definition for a Gaussian TT-tensor is given in~\cite{rakhshan2020tensorized}, but the approach taken here differs considerably in how we use the randomized tensor.

In \Cref{alg:TT-rounding-Rand2}, we first generate a random Gaussian tensor $\T{R}$ with given target TT-ranks $\{\ell_n\}$ following \Cref{def:randtt}. Next, we use the efficient multiplication of tensor $\T{R}$ with a given tensor $\T{Y}$, see \Cref{alg:contractionr2l}, to obtain the sketches (sometimes also referred to as partial random projections) $\{\M{W}_n\}$ of $\T{Y}$ (\emph{randomization phase}). A visualization of this process is provided in \Cref{fig:one-sided}. Finally, we construct a left-orthogonal compressed TT-tensor $\T{X}$. Starting with $n=1$ and $\TT{X}{1} = \TT{Y}{1}$, we compute the QR factorization of the sketched matrix, i.e.,
\[
	\left[\VOp(\TT{X}{n})\HOp(\TT{Y}{n+1:N})\right]\HOp(\TT{R}{n+1:N})^\top = \VOp(\TT{X}{n})\M{W}_n = \M{Q}_n\M{R}_n. 
\]

Since at the $n$th step the first $n-1$ cores of $\T{X}$ are already orthogonalized, they do not need to be considered explicitly in the above factorization. By projecting $\VOp(\TT{X}{n})\HOp(\TT{Y}{n+1:N})$ onto the column space of $\M{Q}_n$, we approximate the product of the final $N-n+1$ cores as
\[
	\VOp(\TT{X}{n})\HOp(\TT{Y}{n+1:N})\approx \M{Q}_n\M{Q}_n^\top \VOp(\TT{X}{n})\HOp(\TT{Y}{n+1:N}) = \M{Q}_n\M{M}_n\HOp(\TT{Y}{n+1:N}).
\]
Then, the cores are updated, i.e., $\VOp(\TT{X}{n}) = \M{Q}_n$ and $\HOp(\TT{X}{n+1}) =  \M{M}_n \HOp(\TT{Y}{n+1})$. 

\smallskip
\noindent
It is important to mention that the Randomize-then-Orthogonalize approach produces a left-orthogonal tensor $\T{X}$. We can use this observation to compress the tensor further. Therefore, if the ranks are not known a priori, then we choose the ranks to be sufficiently large and truncate them further by using \Cref{alg:TT-rounding-RLR}. 
In particular, since the output tensor of \cref{alg:TT-rounding-Rand2} is left orthogonal, we can skip the orthogonalization phase (\cref{line:TTRound:OrthRL} of \Cref{alg:TT-rounding-RLR}), and execute
\crefrange{line:TTRound:Line3}{line:TTRound:Line10}.
This is what we do in our numerical experiments when the rank is not known \textit{a priori}.

\begin{algorithm}[!ht]
\caption{TT-Rounding: Randomize-then-Orthogonalize}
\label{alg:TT-rounding-Rand2}
\begin{algorithmic}[1]
	\Require{A tensor $\T{Y}$ in TT format with ranks $\{R_n\}$, target TT-ranks $\{\ell_n\}$}
	\Ensure{A tensor $\T{X}$ in TT format with ranks $\{\ell_n\}$}
	\Function{$\T{X} =$ TT-Rounding-RandOrth}{$\T{Y}, \{\ell_n\}$}
	    \State Select a random Gaussian TT-tensor $\T{R}$ with target TT-ranks $\{\ell_n\}$
	    \State $\{\M{W}_n\} =$ \Call{PartialContractionsRL}{$\T{Y},\T{R}$} \Comment{compute partial random contractions} 
	    \label{line:prc_rto}
		\State{$\TT{X}{1}=\TT{Y}{1}$}
		\For{$n=1$ to $N-1$}
          \State{$\M{Z}_n = \VOp(\TT{X}{n})$ \Comment{$\TT{X}{n}$ is $\ell_{n-1} \times I_n \times R_n$}} \label{line:unfold_rto}
		      \State{$\M{Y}_n = \M{Z}_n \M{W}_n$ \Comment{form the sketched matrix}}\label{line:sketch_rto}
		    \State{$[\VOp(\TT{X}{n}),\sim] = \textsc{QR}(\M{Y}_n)$ \Comment{thin QR to compute an orthonormal basis}}\label{line:qr_rto}
			\State{$\M{M}_n = \VOp(\TT{X}{n})^\Tra \M{Z}_n$ \Comment{form $\ell_n \times R_n$ matrix}}\label{line:mult1_rto}
			\State{$\HOp(\TT{X}{n+1}) = \M{M}_n \HOp(\TT{Y}{n+1})$ \Comment{ $\TT{X}{n+1} = \TT{Y}{n+1} \times_1 \M{M}_n$}}\label{line:mult2_rto}
		\EndFor
	\EndFunction
\end{algorithmic}
\end{algorithm}

\begin{figure}[!ht]
     \centering
    \scalebox{0.6}{\tikzstyle{vertex} = [circle,draw,scale=.8]
\tikzstyle{edge} = [line width=1pt]

\begin{tikzpicture}[scale=.8]

\node (a0) at (0,0) {};
\node[vertex,fill=blue!40] (a1) at (3,0) {$\hspace{0.08in}\TT{X}{n}\hspace{0.08in}$};
\node[vertex,fill=blue!40] (a2) at (6,0) {$\TT{X}{n+1}$};
\node[vertex,fill=blue!40] (a3) at (9,0) {$\TT{X}{n+2}$};
\node (a4) at (12,0) {$\cdots$};
\node[vertex,fill=blue!40] (a5) at (15,0) {$\TT{X}{N-1}$};
\node[vertex,fill=blue!40] (a6) at (18,0) {$\hspace{0.08in}\TT{X}{N}\hspace{0.08in}$};

\foreach \i in {0,...,6} {
    \node (b\i) at (\i*3,-3) {};
}

\node (c0) at (0,-6) {};
\node (c1) at (4.25,-6) {};
\node[vertex] (c2) at (6,-6) {$\TT{R}{n+1}$};
\node[vertex] (c3) at (9,-6) {$\TT{R}{n+2}$};
\node[gray] (c4) at (12,-6) {$\cdots$};
\node[vertex] (c5) at (15,-6) {$\TT{R}{N-1}$};
\node[vertex] (c6) at (18,-6) {$\hspace{0.08in}\TT{R}{N}\hspace{0.08in}$};

\draw[edge] (a0) -- (a1) node[midway,above] {$R_{n-1}$};
\draw[edge] (a1) -- (a2) node[midway,above] {$R_n$};
\draw[edge] (a2) -- (a3) node[midway,above] {$R_{n+1}$};
\draw[edge] (a3) -- (a4) node[midway,above] {$R_{n+2}$};
\draw[edge] (a4) -- (a5) node[midway,above] {$R_{N-2}$};
\draw[edge] (a5) -- (a6) node[midway,above] {$R_{N-1}$};

\draw[edge] (c1) -- (c2) node[midway,above] {$\ell_n$};
\draw[edge] (c2) -- (c3) node[midway,above] {$\ell_{n+1}$};
\draw[edge] (c3) -- (c4) node[midway,above] {$\ell_{n+2}$};
\draw[edge] (c4) -- (c5) node[midway,above] {$\ell_{N-2}$};
\draw[edge] (c5) -- (c6) node[midway,above] {$\ell_{N-1}$};

\draw[edge] (a1) -- (b1) node[midway,right] {$I_n$};
\draw[edge] (a2) -- (b2) node[midway,right] {$I_{n+1}$};
\draw[edge] (a3) -- (b3) node[midway,right] {$I_{n+2}$};
\draw[edge] (a5) -- (b5) node[midway,right] {$I_{N-1}$};
\draw[edge] (a6) -- (b6) node[midway,right] {$I_N$};

\draw[edge] (c2) -- (b2) node[midway,right] {$I_{n+1}$};
\draw[edge] (c3) -- (b3) node[midway,right] {$I_{n+2}$};
\draw[edge] (c5) -- (b5) node[midway,right] {$I_{N-1}$};
\draw[edge] (c6) -- (b6) node[midway,right] {$I_N$};

\end{tikzpicture}}
     \caption{Random projection for the Randomize-then-Orthogonalize \Cref{alg:TT-rounding-Rand2}.}
     \label{fig:one-sided}
 \end{figure}

\begin{algorithm}[!ht]
\caption{TT-Rounding: Two-Sided-Randomization (Generalized Nystr\"om)}
\label{alg:TT-rounding-Rand3}
\begin{algorithmic}[1]
	\Require{A tensor $\T{Y}$ in TT format with ranks $\{R_n\}$, target TT-ranks $\{\ell_n\}$ and $\{\rho_n\}$}
	\Ensure{A tensor $\T{X}$ in TT format with ranks $\{\ell_n\}$}
	\Function{$\T{X} =$ TT-Rounding-RandOrth}{$\T{Y}, \{\ell_n\}$}
	    \State{Generate random Gaussian TT-tensor $\T{L}$ with ranks $\{\ell_n\}$}
	    \State{Generate random Gaussian TT-tensor $\T{R}$ with ranks $\{\rho_n\}$ \Comment{choose $\rho_n>\ell_n$}}
		\State $\{\M{W}_n^L\} =$ \Call{PartialContractionsLR}{$\T{Y},\T{L}$}
			    \Statex{\Comment{Precompute sketches from the left}} \label{line:tsr_pclr}
	    \State $\{\M{W}_n^R\} =$ \Call{PartialContractionsRL}{$\T{Y},\T{R}$} \label{line:tsr_pcrl}
	     \Statex{\Comment{Precompute sketches from the right}}
		\For{$n=1$ to $N-1$}
		    \State{$[\M{U}_n,\M{\Sigma}_n,\M{V}_n] = \textsc{SVD}(\M{W}_n^L \M{W}_n^R)$ \Comment{Compute SVD of $\M{W}_n^L \M{W}_n^R$}} 
		    \Statex{\Comment{$\M{V}_n$ is $\rho_n \times \ell_n$} \nonumber}\label{line:tsr_svd}
		    \State{$\M{L}_n = \M{W}_n^R \M{V}_n (\M{\Sigma}_n^\Pinv)^{1/2}$ \Comment{Determine internal $R_n \times \ell_n$ left factor $\M{L}_n$}} \label{line:tsr_left}
		    \State{$\M{R}_n = (\M{\Sigma}_n^\Pinv)^{1/2} \M{U}_n^\Tra \M{W}_n^L$ \Comment{Determine internal $\ell_n \times R_n$ right factor $\M{R}_n$}}\label{line:tsr_right}
		\EndFor
		\State{$\VOp(\TT{X}{1}) = \VOp(\TT{Y}{1}) \M{L}_1$}
		\For{$n=2$ to $N-1$} \label{line:tsr_for}
		    \State{$\HOp(\TT{X}{n}) = \M{R}_{n-1} \HOp{\big(\VOp(\TT{Y}{n}) \M{L}_{n}\big)}$ \Comment{$\TT{X}{n} = \TT{Y}{n} \times_1 \M{R}_{n-1} \times_3 \M{L}_n$}}
		    \Statex{\Comment{hence $\TT{X}{n}$ is $\ell_{n-1} \times I_n \times \ell_n$}}
		\EndFor
		\State{$\HOp(\TT{X}{N}) = \M{R}_{N-1} \HOp(\TT{Y}{N})$}
	\EndFunction
\end{algorithmic}
\end{algorithm}

\subsection{Two-Sided-Randomization}\label{ssec:twosided} Analogous to the one-sided \Cref{alg:TT-rounding-Rand2}, we start with generating two TT random Gaussian tensors $\T{L}$ and $\T{R}$ with given target TT-ranks $\{\ell_n\}$ and $\{\rho_n\}$ (with $\rho_n > \ell_n)$ and
computing the sketches  
$\{\M{W}_n^L\}$, $\{\M{W}_n^R\}$ of $\T{Y}$ from the left and right, respectively (\emph{randomization phase}, see \Cref{fig:two-sided}). Next, for each $n = 1, \ldots, N-1$ we compute the SVD of a product of partial contractions $\M{W}_n^L \M{W}_n^{R}$, i.e., $\M{W}_n^L \M{W}_n^{R} = \M{U}_n\M{\Sigma}_n\M{V}_n^\top$,
and form left and right factor matrices
\begin{equation}
\label{eq:LnRn}
\M{L}_n = \M{W}_n^R \M{V}_n (\M{\Sigma}_n^\dagger)^{1/2} \quad \mbox{ and } \quad
\M{R}_n = (\M{\Sigma}_n^\dagger)^{1/2}\M{U}_n^\top\M{W}_n^L.
\end{equation}
In order to highlight the significance of matrices $\M{L}_n$ and $\M{R}_n$, we consider the following unfolding of the TT-tensor $\T{Y}$
\[
	\M{Y}_{(1:n)} = \VOp(\TT{Y}{1:n})\HOp(\TT{Y}{n+1:N}),
\]
with factors $\VOp(\TT{Y}{1:n})\in \R^{(I_1\cdots I_n) \times R_n}$ and $\HOp(\TT{Y}{n+1:N}) \in \R^{R_n \times (I_{n+1}\cdots I_N)}$. Similarly, we define matrices 
\begin{align*}
	\M{\Psi}_n &:= \VOp(\TT{L}{1:n})^\top \in \R^{\ell_n \times (I_1\cdots I_n)}, \\
	\M{\Omega}_n &:= \HOp(\TT{R}{n+1:N})^\top \in \R^{ (I_{n+1}\cdots I_N) \times \rho_n}
\end{align*}
as the partial unfoldings of random Gaussian TT-tensors $\T{L}$ and $\T{R}$, respectively. Then multiplying matrix $\M{Y}_{(1:n)}$ on the left by $\M{\Psi}_n$ and on the right by $\M{\Omega}_n$ yields
\[
\M{\Psi}_n \M{Y}_{(1:n)} \M{\Omega}_n = (\M{\Psi}_n \VOp(\TT{Y}{1:n}))(\HOp(\TT{Y}{n+1:N}) \M{\Omega}_n) = \M{W}_n^L \M{W}_n^R. 
\]

\noindent
Following identity~\eqref{eq:EssenNY} illustrating the main idea of the Generalized Nystr\"{o}m method for matrices discussed in \Cref{ssec:randsvd}, we have
\[
\begin{aligned}
\M{Y}_{(1:n)} & \approx \big(\M{Y}_{(1:n)} \M{\Omega}_n \big) 
\big(\M{\Psi}_n \M{Y}_{(1:n)} \M{\Omega}_n \big)^\dagger \big(\M{\Psi}_n \M{Y}_{(1:n)} \big) \\
& = \VOp(\TT{Y}{1:n}) \M{W}_n^R \big(\M{W}_n^L \M{W}_n^R \big)^\dagger \M{W}_n^L \HOp(\TT{Y}{n+1:N}) \\
&= \VOp(\TT{Y}{1:n}) \M{L}_n\M{R}_n \HOp(\TT{Y}{n+1:N}),
\end{aligned}
\]
see \Cref{fig:two-sided1}.

Having all left and right factors at hand, for each core of the tensor $\T{Y}$ (treating the first and last core separately), we distribute them according
to the formula
\[
\HOp(\TT{X}{n}) = \M{R}_{n-1} \HOp(\VOp(\TT{Y}{n}) \M{L}_{n})
\]
forming the cores of the resulting tensor $\T{X}$, see \Cref{fig:two-sided2}.

In contrast to \Cref{alg:TT-rounding-Rand2},
the Two-Sided-Randomization approach does not produce an orthogonal tensor. However, with a little restructuring, it can be adapted to produce an orthogonal tensor (we do not discuss that here). This variation may be useful for the case when the target TT-ranks are not known in advance, and producing an orthogonal tensor can be used in conjunction with \Cref{alg:TT-rounding-RLR} to further compress the tensor.

\subsection{Rounding of TT-sums}\label{ssec:ttsums}

One of the most common arithmetic operations that depends on TT-rounding is TT-summation, i.e., we want to compress a tensor $\T{Y}$ that is available as the sum of $s$ TT-tensors: $ \T{Y} = \T{Y}^{(1)} + \dots + \T{Y}^{(s)}$.
To reuse existing algorithms, there are two options available to us. For example, we can form the TT-tensor $\T{Y}$ explicitly and then apply one of the compression algorithms proposed previously. As we will argue in \Cref{ssec:compcosts}, the computational costs of this approach has cubic scaling with respect to the number of summands $s$ using the TT-rounding approach, and a quadratic scaling with respect to $s$ using the randomized approaches (Randomize-then-Orthogonalize and Two-Sided-Randomization). This is computationally infeasible as $s$ becomes large. Alternatively, we can form the partial sum $\T{Y}^{(1)} + \T{Y}^{(2)}$, compress this partial sum, add the resulting truncated term to the summand $\T{Y}^{(3)}$ and proceed in the same way with remaining terms. Variations of this approach can be performed using ideas from the summation methods described in~\cite[Chapter 4.1]{higham2002accuracy}. These approaches scale linearly with $s$, but can result in large errors in the compressed representation due to the repeated truncation.

\begin{algorithm}[t]
\caption{TT-Rounding of a Sum: Randomize-then-Orthogonalize}
\label{alg:TT-rounding-sum-Rand2}
\begin{algorithmic}[1]
	\Require{Tensors $\{\T{Y}^{(j)}\}_{1\le j \le s}$ in TT format with ranks $\{R_n^{(j)}\}_{1\leq j\leq s}$, target TT-ranks $\{\ell_n\}$}
	\Ensure{A tensor $\T{X} \approx \sum_{j=1}^s \T{Y}^{(j)}$ in TT format with ranks $\{\ell_n\}$}
	\Function{$\T{X} =$ TT-Rounding-sum-RandOrth}{$\{\T{Y}^{(j)}\}_{1\le j \le s}, \{\ell_n\}$}
	    \State Select random Gaussian TT-tensor $\T{R}$ with ranks $\{\ell_n\}$
	    \For{$j=1$ to $s$}
    	    \State $\left\{\M{W}_n^{(j)}\right\} =$ \Call{PartialContractionsRL}{$\T{Y}^{(j)},\T{R}$} \Comment{precompute sketches from the right} \label{line:sum:contract}
        \EndFor
        \State{$\TT{X}{1} = \begin{bmatrix} \T{T}_{\T{Y}^{(1)},1} & \ldots & \T{T}_{\T{Y}^{(s)},1}\end{bmatrix}$}
		\For{$n=1$ to $N-1$}
		    \State{{$\M{Z}_n = \VOp(\TT{X}{n})$}} \Comment{$\TT{X}{n}$ is $\ell_{n-1} \times I_n \times \sum_{j=1}^s R_{n}^{(j)}$}
		    \State{$\M{Y}_n = \VOp(\TT{X}{n}) \begin{bmatrix} \M{W}_{n}^{(1)} \\ \vdots \\ \M{W}_{n}^{(s)}\end{bmatrix}$ \Comment{complete random sketch}} \label{line:sum:proj}
		    \State{$[\VOp(\TT{X}{n}),\sim] = \textsc{QR}(\M{Y}_n)$ \label{line:sum:qr} \Comment{thin QR factorization}}
			\State{$\begin{bmatrix} \M{M}_{n}^{(1)} & \ldots & \M{M}_{n}^{(s)} \end{bmatrix} = \VOp(\TT{X}{n})^\Tra \M{Z}_n$} \label{line:sum:old}
			\If{$n < N-1$} \Comment{exploit structure in next internal core}
    			\State{$\HOp(\TT{X}{n+1}) = \begin{bmatrix} \M{M}_{n}^{(1)}\HOp\left(\T{T}_{\T{Y}^{(1)},n+1}\right) & \cdots & \M{M}_{n}^{(s)}\HOp\left(\T{T}_{\T{Y}^{(s)},n+1} \right) \end{bmatrix}$} \label{line:sum:next}
            \Else
                \State{$\TT{X}{N} = \displaystyle \sum_{j=1}^s \M{M}_{n-1}^{(j)} \T{T}_{\T{Y}^{(j)},N}$} \label{line:sum:nextlast}
            \EndIf
		\EndFor
	\EndFunction
\end{algorithmic}
\end{algorithm}

In this subsection, we show how to combine the addition and randomized rounding operations to reduce further the computational costs, which is particularly effective when the number of summands $s$ is large. The basic idea is to exploit the nonzero structure of the TT-cores of the sum of TT-tensors to avoid computing with zeros. Applying the orthogonalization phase, as required to perform the deterministic truncation phase, requires assembling the TT representation of the sum. 
Furthermore, the orthogonalization and multiplication by the triangular factor destroys the structure in the middle cores. By using randomization, we can avoid this explicit TT assembly of the sum and avoid unnecessary computations on zeros. \Cref{alg:TT-rounding-sum-Rand2} provides the pseudocode for rounding the sum of $s$ input TT-tensors.
To simplify the notation, we derive the efficient computations considering the case $s=2$, as the generalization will be clear.
Let $\T{Y}$ and $\T{Z}$ be two TT-tensors and consider the TT-tensor $\T{X} = \T{Y} + \T{Z}$.
Let $\T{R}$ be a given random Gaussian TT-tensor and let $\{\M{W}_n^{\T{Y}}\}$ and $\{\M{W}_n^{\T{Z}}\}$, for $n = 1,\ldots,N-1$ be the right-to-left partial contractions of $\T{Y}$ and $\T{Z}$ with $\T{R}$. We have
\[
\M{X}_{(1:n)} = \VOp(\TT{X}{1:n}) \HOp(\TT{X}{n+1:N}) = \begin{bmatrix} \VOp(\TT{Y}{1:n}) & \VOp(\TT{Z}{1:n}) \end{bmatrix} \begin{bmatrix} \HOp(\TT{Y}{n+1:N}) \\ \HOp(\TT{Z}{n+1:N}) \end{bmatrix},
\]
and by using \cref{eq:four_prod_mat_rep},
\begin{align*}
\VOp(\TT{X}{1:n}) &=\begin{bmatrix}(\M{I}_{I_n} \otimes \VOp(\TT{Y}{1:n-1})) & (\M{I}_{I_n} \otimes \VOp(\TT{Z}{1:n-1}))\end{bmatrix} \begin{bmatrix}\VOp(\TT{Y}{n}) & \\ & \VOp(\TT{Z}{n}) \end{bmatrix},\\
\HOp(\TT{X}{n+1:N}) &= \begin{bmatrix}\HOp(\TT{Y}{n+1}) & \\ & \HOp(\TT{Z}{n+1})\end{bmatrix}\begin{bmatrix}\HOp(\TT{Y}{n+2:N}) \otimes \M{I}_{I_{n+1}}\\\HOp(\TT{Z}{n+2:N}) \otimes \M{I}_{I_{n+1}}\end{bmatrix}.
\end{align*}
The matrix $\M{W}_n^{\T{X}}$ can be expressed as 
\[
	\M{W}_n^{\T{X}} = \HOp(\TT{X}{n+1:N})\HOp(\TT{R}{n+1:N})^\top = \begin{bmatrix}
			\HOp(\TT{Y}{n+1:N}) \\ \HOp(\TT{Z}{n+1:N})
		\end{bmatrix}\HOp(\TT{R}{n+1:N})^\top = \begin{bmatrix} \M{W}_n^{\T{Y}} \\ \M{W}_n^{\T{Z}} \end{bmatrix}.
\]
This justifies the procedure in \Cref{alg:TT-rounding-sum-Rand2} of computing the partial contractions separately for each summand (\Cref{line:sum:contract}) and concatenating them (\Cref{line:sum:proj}). 

After the QR factorization of the projected matrix produces the truncated core, we compute the contraction between the new and old cores (\Cref{line:sum:old}) and store the result in a matrix $\M{M}_n = \begin{bmatrix}
	\M{M}_n^{\T{Y}} & \M{M}_n^{\T{Z}}
\end{bmatrix}$ of size $\ell_n \times (R_n^{\T{Y}}+R_n^{\T{Z}})$. 
This matrix is now multiplied from the right by $\HOp(\TT{X}{n+1:N})$ to compute the updated right factor of $\M{X}_{(1:n)}$. This multiplication can be absorbed by the $(n+1)$th core as follows:
\begin{align*}
\M{M}_n \HOp(\TT{X}{n+1:N}) &= \begin{bmatrix}
	\M{M}_n^{\T{Y}} & \M{M}_n^{\T{Z}}\end{bmatrix}  \begin{bmatrix}\HOp(\TT{Y}{n+1}) & \\ & \HOp(\TT{Z}{n+1})\end{bmatrix}\begin{bmatrix}\HOp(\TT{Y}{n+2:N}) \otimes \M{I}_{I_{n+1}}\\\HOp(\TT{Z}{n+2:N}) \otimes \M{I}_{I_{n+1}}\end{bmatrix},\\
	&= \begin{bmatrix}\M{M}_n^{\T{Y}}\HOp(\TT{Y}{n+1})  & \M{M}_n^{\T{Z}}\HOp(\TT{Z}{n+1})\end{bmatrix}\begin{bmatrix}\HOp(\TT{Y}{n+2:N}) \otimes \M{I}_{I_{n+1}}\\\HOp(\TT{Z}{n+2:N}) \otimes \M{I}_{I_{n+1}}\end{bmatrix}.\\
\end{align*}
Hence we update $\HOp{\left(\TT{X}{n+1}\right)}$ as in \Cref{line:sum:next}.
For $n=N-1$ we have
\begin{equation*}
    \TT{X}{N} = \begin{bmatrix} \M{M}_{N-1}^{\T{y}} & \M{M}_{N-1}^{\T{z}}\end{bmatrix} \begin{bmatrix} \TT{Y}{N} \\ \TT{Z}{N} \end{bmatrix}\\
              = \M{M}_{N-1}^{\T{y}} \TT{Y}{N} + \M{M}_{N-1}^{\T{z}} \TT{Z}{N}.
\end{equation*}
Generalizing this expression to $s$ terms yields \cref{line:sum:nextlast}.

\subsection{Computational costs} \label{ssec:compcosts}
To analyze the computational cost, we make the following assumptions that will simplify the analysis. Let $\T{Y} \in \mathbb{R}^{I \times \dots \times I}$ be a tensor of order $N$  with ranks $(1,R, \dots, R,1)$ in TT format. We want to compress $\T{Y}$ to obtain a TT-tensor $\T{X}$ with ranks $(1,\ell, \dots, \ell,1)$. 
Here and in \Cref{sec:CompCostStandard}, we assume that $\ell = \Theta(R)$.

\subsubsection{Randomized compression algorithms}\label{subsec: RCA}
We now analyze the computational cost of \Cref{alg:TT-rounding-Rand1,}, \Cref{alg:TT-rounding-Rand2}, and \Cref{alg:TT-rounding-Rand3}. 
\paragraph{Orthogonalize-then-Randomize (\Cref{alg:TT-rounding-Rand1})}
We denote the overall computational cost of \Cref{alg:TT-rounding-Rand1} by $C_{\rm OtR}$. 
\cref{line:Rand1:orth} of \Cref{alg:TT-rounding-Rand1} invokes \Cref{alg:orthogonalizationr2l}, and its cost is $C_{\rm Orth}$, see \Cref{eq:cost_R2LO}. 
We now analyze the cost of the loop in \crefrange{line:Rand1:for}{line:Rand1:endfor}.
In \cref{line:Rand1:sketch}, we multiply the matrix $\VOp(\TT{X}{n})$ of size $I\ell \times R$ with the random matrix of size $R\times \ell$; this costs $2IR\ell^2 $ flops.  
In \cref{line:Rand1:qr}, we compute the thin QR factorization of the matrix $\M{Y}_n$ of size $I\ell \times \ell$, which costs $4I\ell^3 + \mathcal{O}(\ell^3)$ flops. 
Finally, in \cref{line:Rand1:proj1,line:Rand1:proj2}, we have two different matrix-matrix multiplications which cost $2IR \ell^2 + 2IR^2\ell$ flops, respectively. The total computational cost of \Cref{{alg:TT-rounding-Rand1}} is 
\begin{eqnarray*}\label{eq: cost_TT-rounding-Rand1}
C_{\rm OtR} = & \> C_{\rm Orth}  + (N-2)(2IR^2\ell + 4IR\ell^2 + 4I \ell^3) + \mathcal{O}(IR^2 + NR^3) \\
= & \> I(N-2) \cdot( 5R^3 + 2R^2 \ell + 4R\ell^2 + 4\ell^3) + \mathcal{O}(IR^2 + NR^3) \quad \text{flops}.
\end{eqnarray*}

 \paragraph{Randomize-then-Orthogonalize (\Cref{alg:TT-rounding-Rand2})}  Denoting the total cost of \Cref{alg:TT-rounding-Rand2} with $C_{\rm RtO}$, we analyze the main components that contribute to the total computational cost. 
 \Cref{line:prc_rto} invokes \Cref{alg:contractionr2l} with the corresponding cost denoted by $C_{\rm Contr}$, see \Cref{eq:cost_R2LC}.
 \Cref{line:sketch_rto,line:mult1_rto} contribute by a factor of $2IR\ell^2$ that is the cost of performing the multiplication $\VOp(\TT{X}{n})\M{W}_n$ of sizes $\ell I \times R$ with $R \times \ell$ and $\VOp(\TT{X}{n})^T \M{Z}_n$ of sizes $\ell \times I \ell$ and $I\ell \times R$ that prepare the matrices $\M{Y}_n$ and $\M{M}_n$, respectively, for the next steps.  
 The term $4I\ell^3 + \mathcal{O}(\ell^3)$ represents the cost of the thin QR factorization, in \Cref{line:qr_rto}, of the matrix $\M{Y}_n$ of size $I\ell \times \ell$.
\Cref{line:mult2_rto} that involves multiplication of matrices of size $\ell \times R$ and $R \times IR$ which costs $2IR^2\ell$ flops. 
 The total cost of \Cref{alg:TT-rounding-Rand2} is 
 \begin{eqnarray*}\label{eq: cost_TT-rounding-Rand2}
C_{\rm RtO} = &\> C_{\rm Contr} + (N-2)(2IR^2\ell + 4 I R \ell^2 + 4I \ell^3) + \mathcal{O}( I R^2 + NR^3) \\
=& \> I (N-2)\cdot(4R^2\ell + 6R \ell^2 + 4\ell^3) + \mathcal{O}( I R^2 + NR^3)\quad \text{flops}.
 \end{eqnarray*}
\paragraph{Two-Sided-Randomization (\Cref{alg:TT-rounding-Rand3})} Here we analyze the total cost of \Cref{alg:TT-rounding-Rand3}  and denote it by $C_{\rm 2SR}$. The cost of \Cref{line:tsr_pclr,line:tsr_pcrl} is 2$C_{\rm Contr}$. Here, we recall that the cost of the partial contractions from the left and from the right is the same since we assume that ranks are the same, i.e., $\rho_n = \ell_n = \ell$. However,  in practice, we choose $\rho_n  > \ell_n$ to be different for numerical stability.
There are two other contributions to the cost that  come from two different matrix multiplication: first, of sizes $\ell \times R$ and $R \times IR$, and second, of sizes $\ell \times R$ and $R \times I \ell$. The cost of the for loop starting with \Cref{line:tsr_for} is $\mathcal{O}(NR^2)$. 
Hence, the total computational cost of \Cref{alg:TT-rounding-Rand3} is 
\begin{eqnarray*} C_{\rm 2SR} = & \> 2\cdot C_{\rm Contr} + (N-2)(2IR^2\ell + 2IR\ell^2) + \mathcal{O}(NR^2) \\
= & \> I(N-2) \cdot(6R^2\ell + 6R\ell^2) + \mathcal{O}(IR\ell + NR^2)\quad \text{flops}.
\end{eqnarray*}
Nakatsukasa~\cite{nakatsukasa2020fast} suggests oversampling from the left, i.e., taking $\rho_n = \lceil 1.5 \ell_n \rceil$. We follow this suggestion in our numerical experiments. With this assumption, the cost is slightly higher, i.e., $I(N-2) \cdot (7R^2\ell + 8.5R\ell^2) + \mathcal{O}(IR\ell + NR^2)$ flops, due to the increased cost of the contraction (\Cref{alg:contractionr2l}) with a larger random tensor. 

\begin{table}[h!]
\caption{Summary of the computational costs (discarding lower order terms) of the randomized algorithms proposed in this paper. For completeness, we also include the computational costs of the deterministic algorithms in \Cref{ssec:ttrounding}. Orth and Contr refer to \Cref{alg:orthogonalizationr2l} and \Cref{alg:contractionr2l}, respectively. }
	\label{tab:compcosts}
	\centering
	\scalebox{0.89}{
	\begin{tabular}{c|l|l}
	Algorithms	& Computational cost (flops) & Simplified Cost (flops)\\ \hline\hline
	Orth & $(N-2)I(5R^3)$ &  $-$ \\
	Contr & $(N-2)I(2R^2 \ell + 2 R\ell^2)$ & $-$ \\
	TT-Rounding & $(N-2)I(5 R^3 +  6R^2\ell + 2R\ell^2 ) $ & $(N-2)IR^3 (5+ 6\beta + 2\beta^2)$ \\ \hline
Orth-then-Rand   &  $(N-2)I(5R^3 + 2R^2\ell + 4R\ell^2 + 4 \ell^3) $ & $(N-2)IR^3 (5 + 2\beta + 4\beta^2 + 4\beta^3)$ \\
	Rand-then-Orth & $ (N-2)I(2R^2\ell + 4 R \ell^2 + 4\ell^3)$ & $(N-2)IR^3 (4\beta + 6\beta^2 + 4\beta^3)$ \\
	Two-Sided-Rand & $(N-2)I(6R^2\ell + 6R\ell^2)$ & $(N-2)IR^3 ( 6\beta + 6\beta^2)$
	\end{tabular}
	}
\end{table} 

\begin{figure}[!ht]
    \centering
    \includegraphics[scale=0.35]{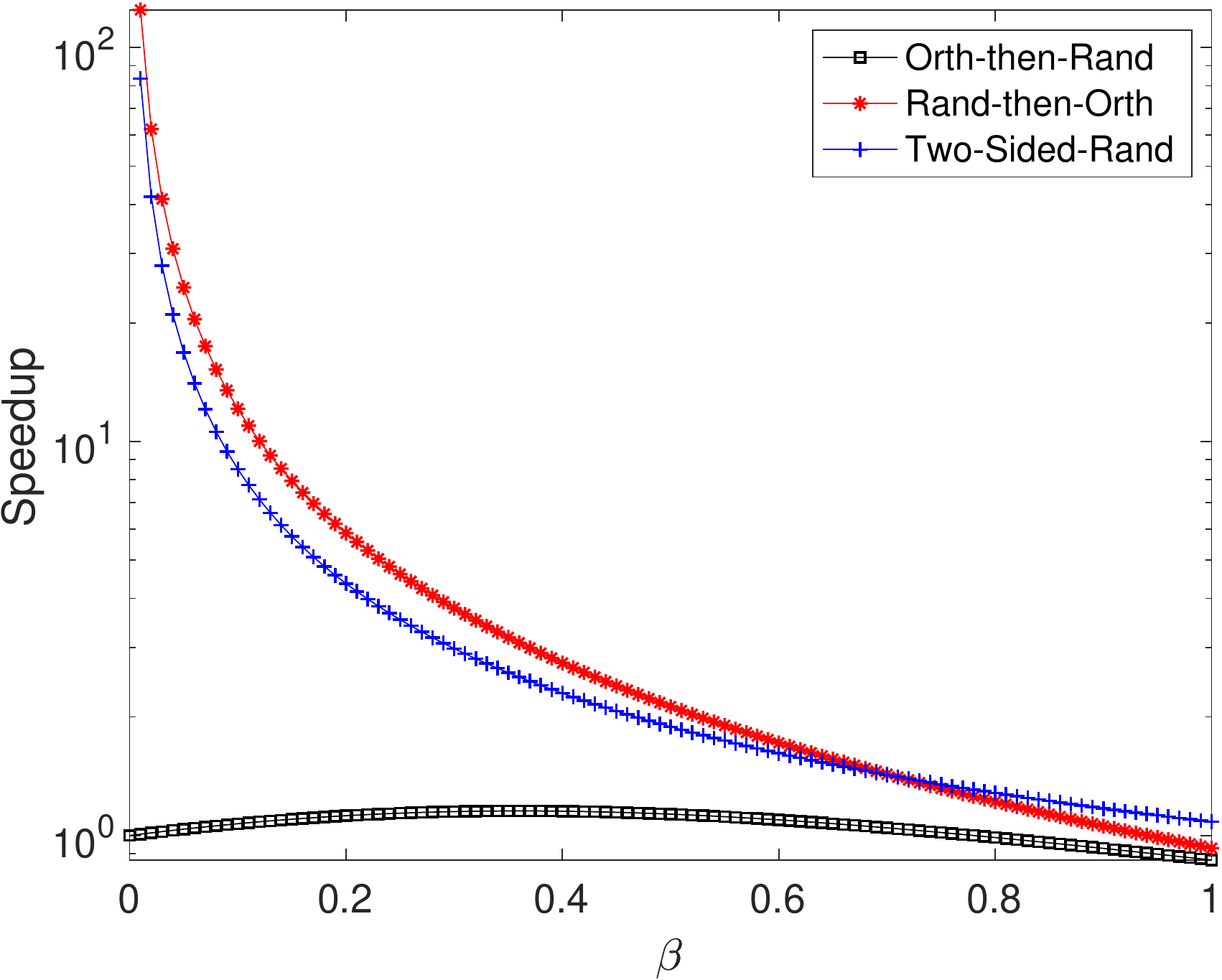}
    \caption{Illustration of the speedups obtained by the randomized algorithms compared with the TT-rounding. Here $\beta = \ell/R$ is the ratio between the target rank and the {original} rank of the tensor. The speedup computations are based on the simplified cost analysis in \Cref{tab:compcosts}.}     \label{fig:compcosts}
\end{figure}

\paragraph{Comparison of different algorithms} To enable the comparison of the different algorithms, we set a target rank  
as $\ell = \beta R$, where $\beta \in (0,1]$ is the ratio between the target rank $\ell$ and the current rank $R$. This allows us to compare the different algorithms more clearly. A summary of the dominant costs of the algorithms is provided in \Cref{tab:compcosts}. For simplicity, we also provide a simplified representation of the computational costs with $\beta = \ell/R$. In \Cref{fig:compcosts}, we plot the speedup of the randomized algorithms compared to the TT-Rounding algorithm; we used the simplified representation of the costs while generating the figure. It can be easily seen that all the proposed methods are faster than the TT-Rounding algorithm. However, the speedup using the Orthogonalize-then-Randomize algorithm is very incremental. In contrast, the other two algorithms, Randomize-then-Orthogonalize and Two-Sided-Randomization, have very similar costs (Randomize-then-Orthogonalize is slightly more efficient for smaller $\beta$) and have much higher speedups especially if $\beta \ll 1$. If $\beta \approx 1 $, both algorithms are very close to the TT-Rounding algorithm. Therefore, the proposed methods are most efficient if $\beta \ll 1$, i.e., the target rank $\ell$ is much smaller compared to the {original} rank $R$. 
\subsubsection{Rounding of TT-sums}To explain the benefits of the algorithm for rounding TT-sums in \Cref{ssec:ttsums}, consider the summation of $s$ tensors of order $N$ (size $I$ in each dimension) each with TT-ranks $(1,R,\dots,R,1)$. Suppose we form the TT-tensor $\T{Y}$, which represents the summation $\T{Y} = \sum_{j=1}^s\T{Y}^{(j)}$, explicitly. The intermediate cores have size $sR \times I \times sR$; the first core is of size $I \times sR$ and the last core is of size $sR\times I$. Suppose the target compression rank in each mode is $\ell$. To leading order, the cost of executing TT-Rounding and Orthogonalize-then-Randomize is $\mathcal{O}(Ns^3R^3I)$. In contrast, the cost of using Randomize-then-Orthogonalize and the Two-Sided-Randomization approach are both $\mathcal{O}(N \ell s^2R^2I)$. This can be beneficial if the number of summands $s$ is large, or the rank $R$ is large. This simple cost analysis does not take into account any structure present in the summation.

Carefully exploiting the structure, as in \Cref{ssec:ttsums}, can reduce this cost.
In particular, by using TT-Rounding of a sum with $s$ tensors of order $N$ with Randomize-then-Orthogonalize summarized in \Cref{alg:TT-rounding-sum-Rand2} the leading order computational cost is $\mathcal{O}(N\ell s R^2 I)$ flops. Notice that by exploiting the structure of the tensor and using randomization, the leading order of the cost is decreased to be linear in $s$ in contrast to cubic in $s$ when no structure was taken into account and the TT-sum tensor was formed explicitly. This decrease in the computational cost is obviously more pronounced when the number of summands $s$ is large. In what follows, we present the analysis of the computational cost of computing the sum of $s$ TT-tensors by randomization and by exploiting the underlying tensor structure.

\paragraph{TT-rounding of a sum: Randomize-then-Orthogonalize (\Cref{alg:TT-rounding-sum-Rand2})}
We analyze the computational cost of \Cref{alg:TT-rounding-sum-Rand2} which we denote by $C_{\rm RtOsum}$. The leading order term is sourced from two main contributions (1) \cref{line:sum:contract} that executes \Cref{alg:contractionr2l} $s$ times resulting in a total computational cost of $sI(N-2)(2R\ell^2+2R^2\ell)$ flops and (2) \cref{line:sum:next} that represents $s$ multiplications between matrices of size $IR\times R$ and $R \times \ell$ resulting in a total cost of $2sIR^2\ell$ flops.
Next we analyze the source of the second leading order term present in the total cost. 
\Cref{line:sum:proj} contributes by a factor of two to the second leading term with a cost $2sIR\ell^2$ flops resulting from multiplying matrices of size $I\ell \times s R$ and $sR\times\ell$. 
Another factor of two comes from the multiplication of two matrices of size $\ell \times I \ell$ and $I \ell \times s R$ in \cref{line:sum:old}, resulting in a total cost of $2sIR\ell^2$ flops. The last factor of two comes from the second leading order term of \Cref{alg:contractionr2l}. Computing the thin QR factorization of the matrix of size $I \ell \times \ell$ in \cref{line:sum:qr} costs $4I\ell^3$ flops. Hence, the total cost of \Cref{alg:TT-rounding-sum-Rand2} is 
\begin{eqnarray*}\label{eq: cost_TT-rounding-sum-Rand2}
C_{\rm RtOsum} = & \> s\cdot C_{\rm Contr} + (N-2)(2sIR^2\ell + 4sIR\ell^2 + 4I\ell^3) + \mathcal{O}(NIR^2) \\
= & \> I(N-2) \cdot(4sR^2\ell + 6sR\ell
^2 + 4\ell^2) + \mathcal{O}(NIR^2)\quad \text{flops}. 
\end{eqnarray*}

\section{Numerical results}\label{sec:numex}
In this section, we illustrate numerically the performance of the newly developed algorithms using tensor data in TT format. We consider both synthetic as well as more realistic test examples. Additional numerical experiments are available in \Cref{sec:hilbert}. All the numerical experiments were performed on {\sc Matlab} R2021a running on a laptop computer with CPU Intel(R) Core(TM) i9-9980H and $64$GB of RAM, using multithreading with $4$ computational threads.

\subsection{TT-tensor with a fixed target rank}

\begin{figure}[ht!]
    \centering
    \includegraphics[width=\textwidth]{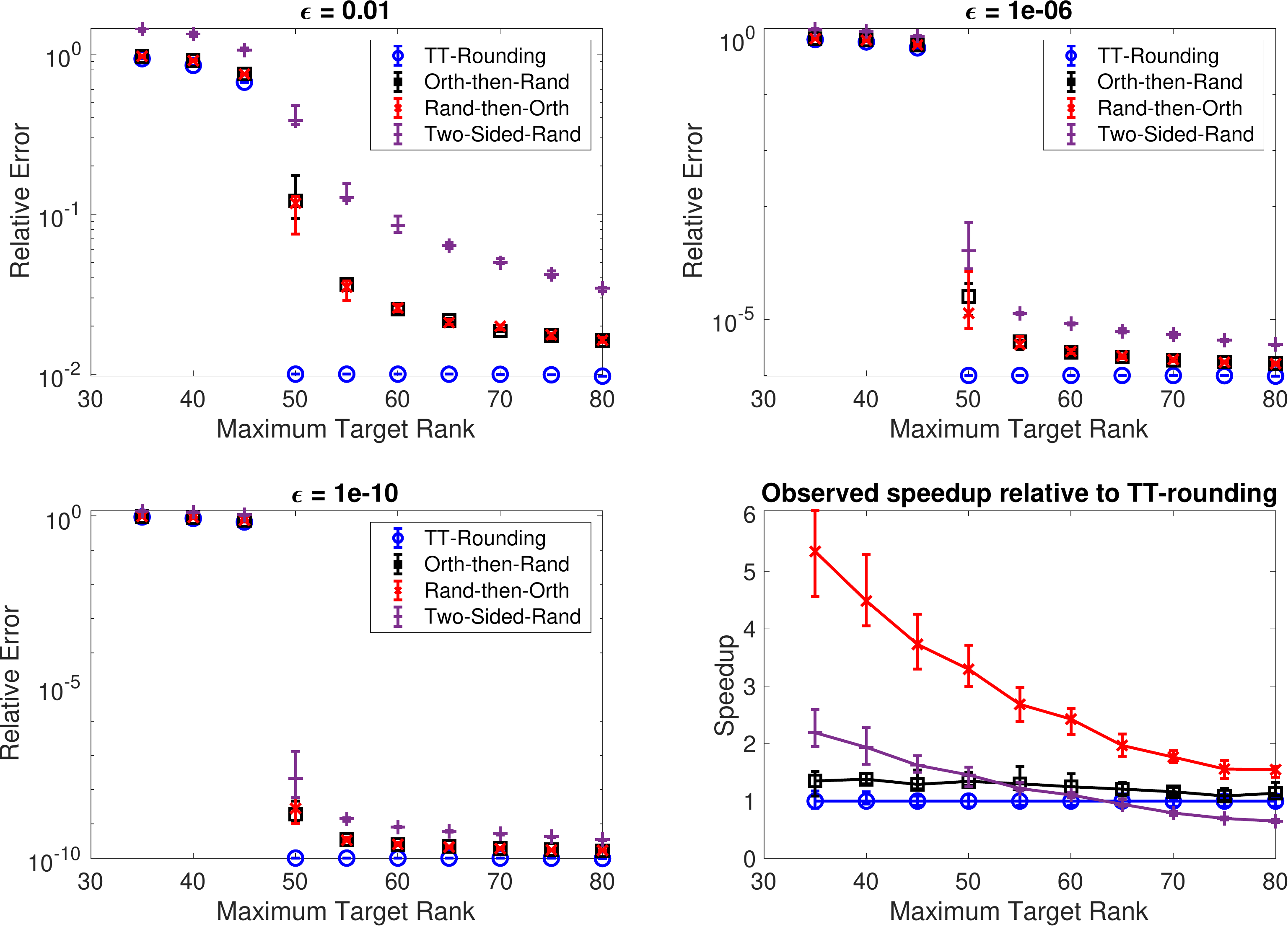}
    \caption{Comparison of error between a low-rank tensor and full rank perturbed tensor, and timings using the deterministic and randomized TT-rounding algorithms for $\epsilon$ perturbed tensor for different values of perturbation $\epsilon$. Statistics were based on $5$ independent runs.}
    \label{Fig:VaryingEpsilon}
\end{figure}

{In our first numerical experiment, we illustrate the performance of our rounding algorithms by rounding a random TT-tensor with a known low-rank representation. 
Throughout, we choose the ranks of the right side randomization in the Two-Sided-Randomization approach (\Cref{alg:TT-rounding-Rand3}) to be $\rho = \lceil 1.5\ell \rceil$ as discussed in \Cref{subsec: RCA}.}

{The random TT-tensor $\T{X}$ is constructed by perturbing a random TT-tensor $\T{X}_1$ with the random TT-tensor $\epsilon \T{X}_2$ as follows: 
$
    \T{X}=\T{X}_1+\epsilon \T{X}_2.
$
TT-tensors $\T{X}_1, \ \T{X}_2\in \mathbb{R}^{100\times\cdots\times100}$ are order $N=10$ normalized random TT-tensors of ranks $(1,50,...,50,1)$ (normalized according to their dimension as described in Definition 3.1), and $\epsilon$ is a perturbation scalar taking the values $\epsilon \in \{10^{-2},10^{-6},10^{-10}\}$. The  ranks  of  the  perturbed  tensor $\T{X}$ are $(1,50 + 50,...,50 + 50,1)$, and the perturbation parameter $\epsilon$ determines how well tensor $\T{X}$ is approximated by the lower rank tensor $\T{X}_1$, i.e., if $\epsilon$ is small, then $\T{X}_1$ is a good rank-$(1,50,...,50,1)$ approximation of  $\T{X}$.}

{We  round the random TT-tensor $\T{X}$ using \Cref{alg:TT-rounding-RLR,alg:TT-rounding-Rand1,alg:TT-rounding-Rand2,alg:TT-rounding-Rand3} to have ranks $(1,\ell,\dots,\ell,1)$, where we vary the parameter $\ell$ from $35$ to $80$ by an increment of $5$. We present these results in \Cref{Fig:VaryingEpsilon}. The approximation error is the relative norm error between tensor $\T{X}$ and the approximate rounded tensor $\widehat{\T{X}}$, i.e., $\frac{\|\T{X} - \widehat{\T{X}}\|}{\|\T{X}\|}$. We also present the time speedup (computed with the average of $5$ independent runs) of the randomized algorithms compared to the deterministic algorithm.}

{For all values of perturbation $\epsilon$, the error resulting from the deterministic algorithm (\Cref{alg:TT-rounding-RLR}) decreases slightly until the ranks of the rounded tensor are $\ell=50$. When $\ell > 50$, the error appears to be very close to $\epsilon$. The errors resulting from the randomized algorithms (\Cref{alg:TT-rounding-Rand1,alg:TT-rounding-Rand2,alg:TT-rounding-Rand3}) are greater than the error resulting from the deterministic algorithm. Additionally, the randomized algorithms produce a more accurate approximation when the target rank is larger and are more accurate when $\epsilon$ is small. The Orthogonalize-then-Randomize and the  Randomize-then-Orthogonalize method (\Cref{alg:TT-rounding-Rand1,alg:TT-rounding-Rand2}) produce a rounded tensor with similar levels of accuracy while the Two-Sided-Randomization approach is the least accurate.  For smaller values of $\epsilon$, there is less difference in the accuracy between the different algorithms. The Randomize-then-Orthogonalize algorithm is the fastest compared to the deterministic algorithm followed by the Two-Sided-Randomization and Orthogonalize-then-Randomize algorithms.}

\subsection{Solving a parametric PDE in the TT format}
\label{sec:PDE_Example_experiments}

\noindent
As a realistic test example, we consider the parameter dependent PDE referred in the literature as the \emph{cookie problem} \cite{KreT11,Tob12}, see \Cref{sec:PDE_Example} for details. Since it is known that the set of solutions of problem~\eqref{eq:cookie} admits a low-rank representation \cite{Gra04,DahDGS16}, we consider a global linear system encapsulating all these linear systems, i.e., 
\[
\left(\sum_{i=1}^N \M{A}_{i,1} \otimes \ldots \otimes \M{A}_{i,N}\right) \T{x} = \T{f},
\]
where $\M{A}_{1,1}$ is the discretization of the operator over the spatial domain with constant parameter values, for each $2\leq i\leq N$, $\M{A}_{i,1}$ is the discretization of the operator over the domain multiplied by the characteristic function corresponding to the corresponding subdomain and $\M{A}_{i,i}$ is a diagonal matrix containing the parameter values for the corresponding parameter, and for each $2 \le j \neq i \le N$, $\M{A}_{i,j}$ is the identity matrix. 
We use the TT-GMRES algorithm \cite{Dol13} to solve this global linear system of equations.
The preconditioned TT-GMRES algorithm builds the basis vectors in TT format $\T{V}_1, \T{V}_2, \dots$ using the inexact Arnoldi procedure; since at each step, the corresponding TT-tensors are rounded, this results in an inexact Krylov subspace method. 
The main bottleneck is the computation of two linear combinations in each iteration. First, the following sum of $N$ tensors with the same ranks as $\T{V}_k$ is formed when applying the operator to the $k$-th basis vector computed at the previous iteration, i.e.,
\[
   \T{W} = \sum_{i=1}^N \left(\M{A}_{i,1} \otimes \ldots \otimes \M{A}_{i,N}\right) \T{Y},
\] 
after application of the preconditioner, $\T{Y} = \left ( \left ( \sum_{i=1}^N \M{A}_{i,1} \right ) ^{-1} \otimes I \otimes \dots \otimes I \right ) \T{V_k}$.
Second, a linear combination of $k+1$ tensors appears when using the Gram-Schmidt algorithm to orthogonalize $\T{W}$ with respect to the previous basis vectors, 
\[
   \T{Z} = \left ( \T{W} - \sum_{j=1}^k h_{jk} \T{V}_k \right ), \qquad \T{V}_{k+1} = \frac{1}{h_{k+1,k}} \T{Z}, \quad \begin{cases}
      h_{jk} = \langle \T{V}_{k}, \T{W} \rangle, & j = 1, \dots, k, \\
      h_{k+1,k} = \Vert \T{Z} \Vert_F.
   \end{cases}
\]
In both cases, the addition of TT-tensors is followed by a TT-rounding operation in order to reduce the ranks and keep the computations tractable. 
Hence, these steps are amenable to acceleration by using the randomized \Cref{alg:TT-rounding-sum-Rand2} as a rounding procedure in the aforementioned computations. 
Because the second linear combination involves $k+1$ summands, the randomized implementation reduces greatly the cost of later TT-GMRES iterations in particular, as the leading order of its computational cost is decreased from cubic to linear in $k$, as detailed in \Cref{ssec:compcosts}.

\begin{figure}[t]
    \centering
    \includegraphics[width=.45\textwidth]{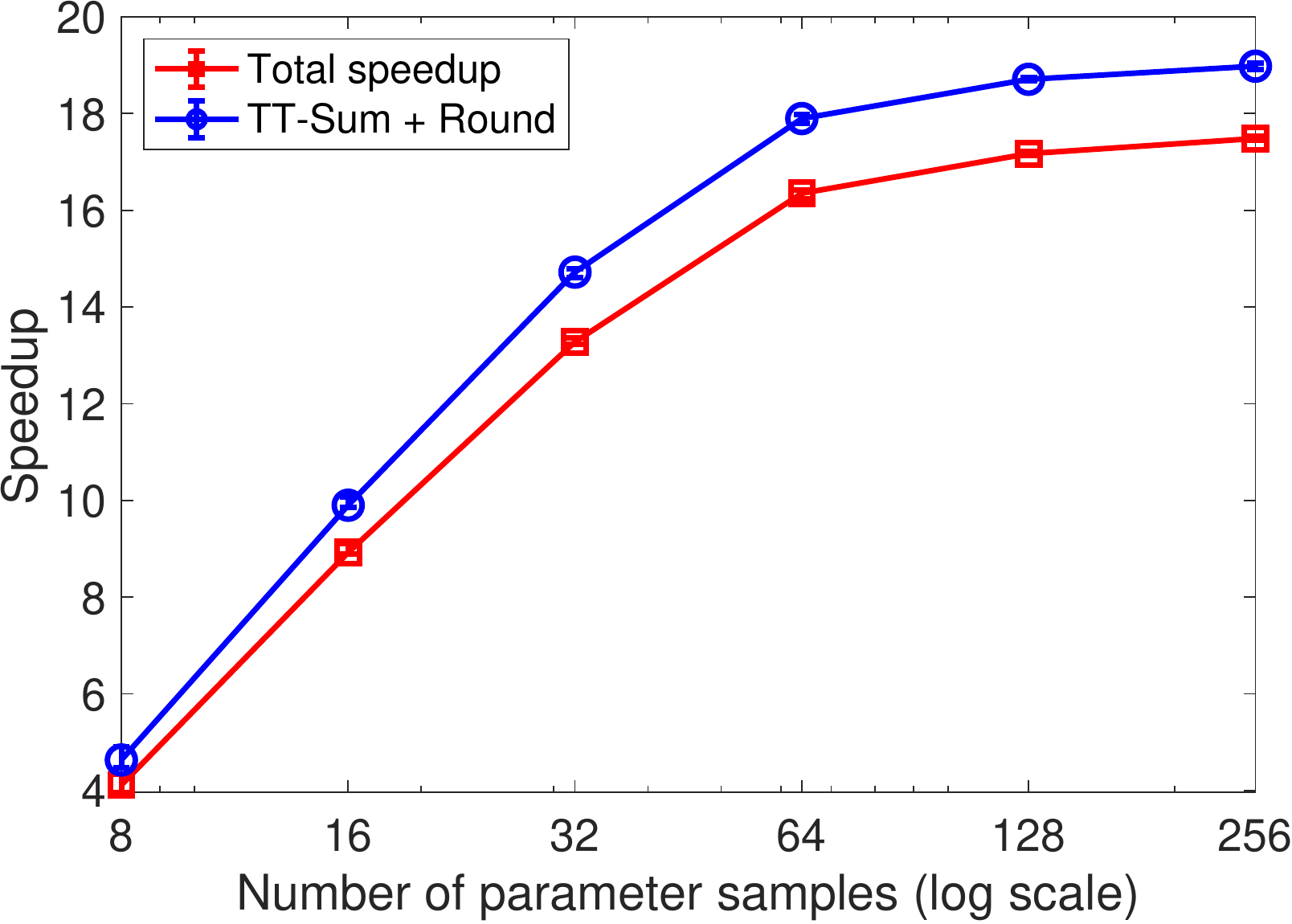} \hspace{.02\textwidth}
    \includegraphics[width=.45\textwidth]{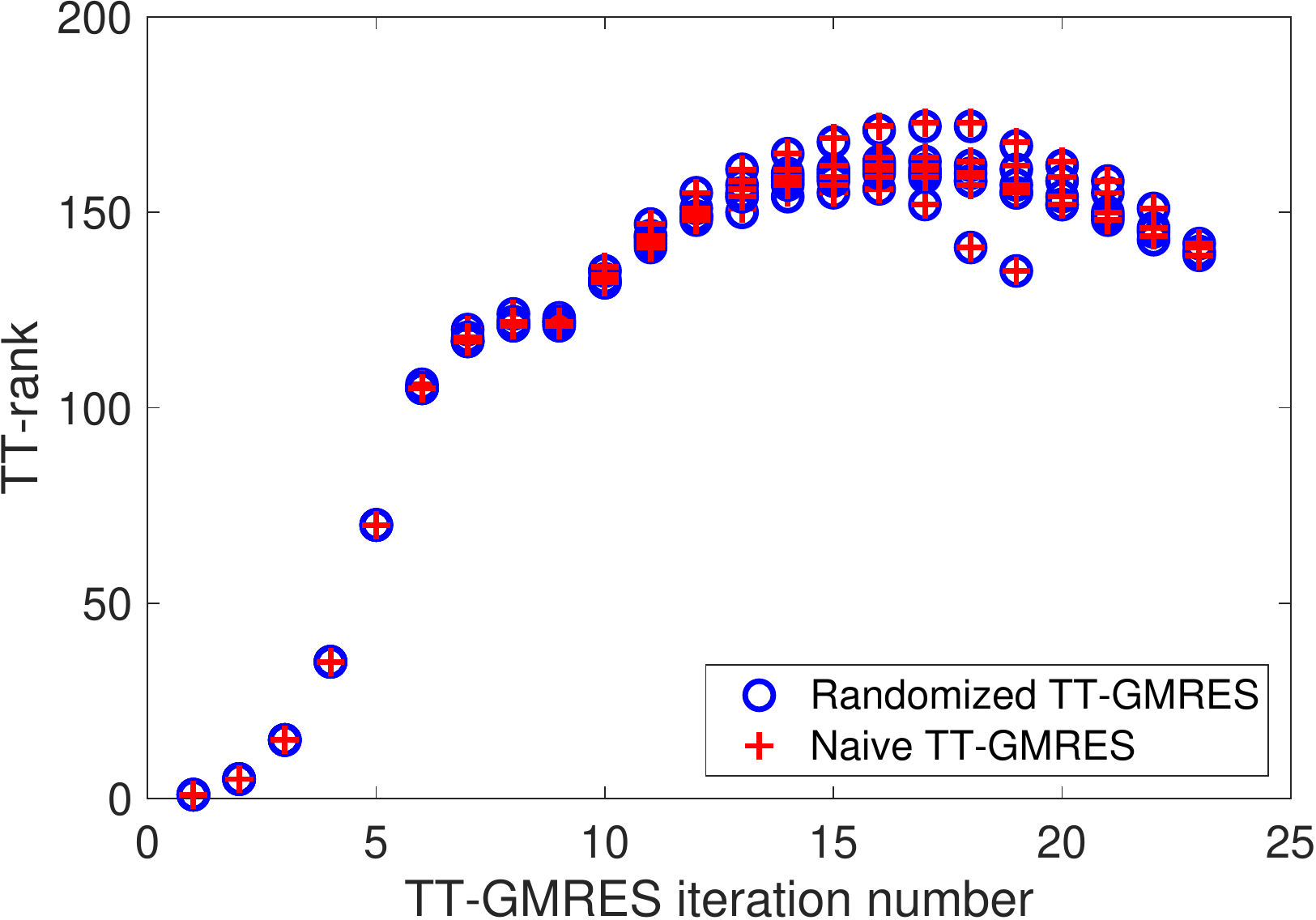}
    \caption{Illustration of the speedups (left) and TT-ranks of the Krylov basis vectors (right) obtained by the deterministic and randomized {summation and rounding} algorithms within the TT-GMRES algorithm to solve problem~\eqref{eq:cookie}.}
    \label{Fig:TTGMRESspeedup}
\end{figure}

\begin{figure}[h!]
    \centering
    \includegraphics[width=0.9\textwidth]{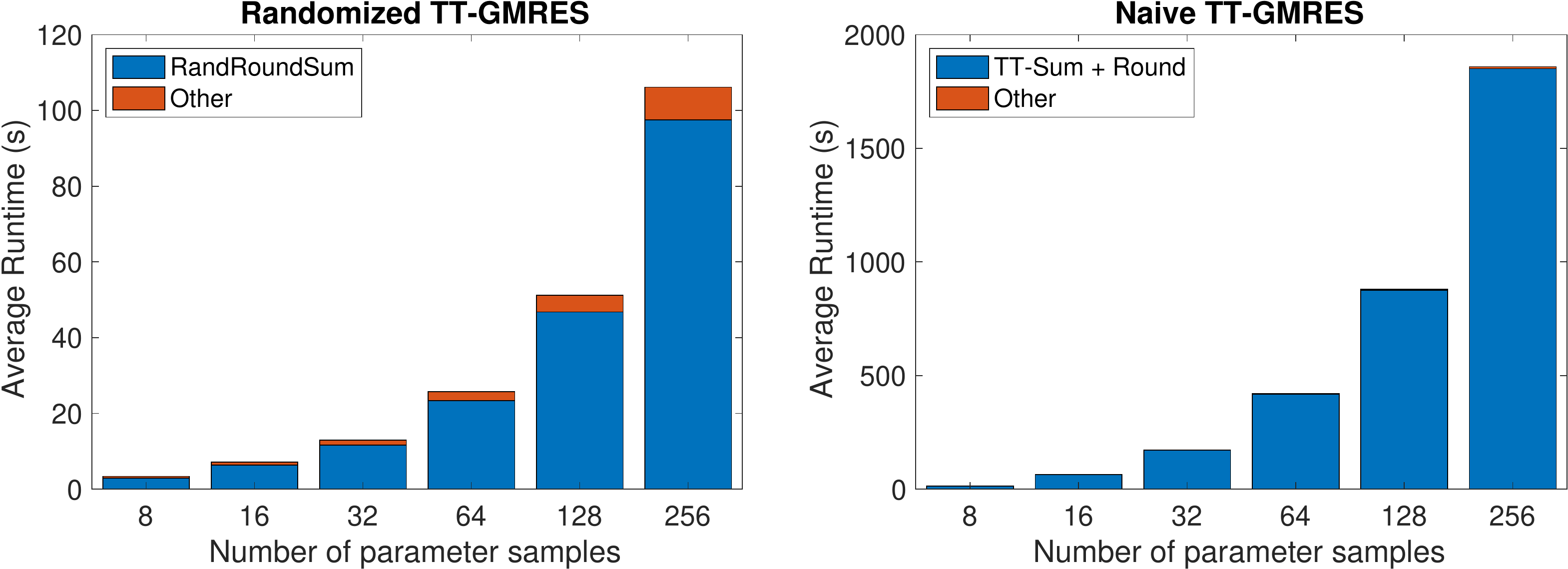}
    \caption{Illustration of the timings using the deterministic and randomized {summation and rounding} algorithms within the TT-GMRES algorithm to solve problem~\eqref{eq:cookie}.}
    \label{Fig:TTGMRESruntime}
\end{figure}

We perform numerical experiments with $N=5$ using a piecewise linear finite element discretization with the mesh presented in \Cref{Fig:cookies}, for various choices of the number of parameter samples, $I = I_2 = \ldots = I_N$, with values of $\rho_i$ distributed linearly between $1$ and $10$. 
The relative tolerance of the TT-GMRES solver is set to $10^{-8}$.
We compare the naive, deterministic implementation of the preconditioned TT-GMRES algorithm with one using randomized summation and rounding steps. Results of the comparison are reported in \Cref{Fig:TTGMRESspeedup}. 
On the right, we display all internal ranks of the TT representation of the Krylov vector computed at that iteration.
We observe that the ranks of the basis vectors and number of iterations are the same using both implementations, and they depend only weakly on the dimension $I$ of the parameter modes of the tensors. 
The speedup achieved by the randomized approach increases consistently with the number of parameter samples.
Taking a closer look at the timing statistics as presented in \Cref{Fig:TTGMRESruntime}, we note that the computation and rounding of the linear combinations identified above indeed dominates the computational cost in both cases, and is a clear computational bottleneck for the deterministic implementation in particular, as the ranks of the sum tensor increase to as much as 2491 in these experiments. 
This explains the remarkable speedup obtained with the randomized approach.

\section{Conclusions and outlook}\label{sec: conclusions}

In this paper, we present randomized algorithms for rounding a tensor, assuming that we have an initial representation in the TT format. This initial representation may not be optimal in terms of storage, and the randomized compression techniques can be used to obtain a more efficient representation. We derive three different algorithms: Orthogonalize-then-Randomize, Randomize-then-Orthogonalize, and Two-Sided-Randomization. We study the computational cost of these algorithms in some detail and show that it can be much smaller than the standard TT-rounding algorithm. Additionally, we consider the special case of rounding a TT-tensor that is represented as the sum of many TT-tensors. 

While applying each of the randomized algorithms proposed here can reduce the computational cost over standard TT-rounding, we further exploit the structure of the problem to reduce the computational cost to be linear in the number of summands. We perform extensive numerical experiments and achieve over 20$\times$ speedups on test problems compared to standard algorithms. There are many avenues for future investigations. First, it would be interesting to derive probabilistic bounds for the accuracy of the rounding approach.  Second, we could consider extending our algorithms to the case where the TT-tensor is obtained as the Hadamard (or elementwise) product of two tensors. Finally, another extension worth considering is developing randomized rounding algorithms in the $\mathcal{H}-$Tucker format.

\section*{Acknowledgements}

This work was initiated as a part of the Statistical and
Applied Mathematical Sciences Institute (SAMSI) Program on Numerical Analysis in Data Science in Fall 2020. Any opinions, findings, and conclusions or recommendations expressed
in this material are those of the authors and do not necessarily reflect the views of the National Science Foundation (NSF). 
Additionally, the authors would like to acknowledge partial support through the NSF: G.B.\ (CCF-1942892), P.C.\ (DMS-1819220), E.H.\ (DMS-1745654), A.M.\ (CCF-1812927), M.P.\ (DMS-1502640), T.W.R.\ (DMS-1745654), A.K.S.\ (DMS-1821149). The authors thank Dr.\ Jocelyn Chi and Prof.\ Eric de Sturler  for the constructive discussions and suggestions.

\newpage
\appendix


\section{Computational cost of standard TT algorithms}
\label{sec:CompCostStandard}

In this section, we analyze the computational cost of the standard approach to TT-rounding which is presented in \Cref{ssec:ttrounding}.

To analyze the computational cost, we make the following assumptions that will simplify the analysis. Let $\T{Y} \in \mathbb{R}^{I \times \dots \times I}$ be a tensor of order $N$  with ranks $(1,R, \dots, R,1)$ in TT format. We want to compress $\T{Y}$ to obtain a TT-tensor $\T{X}$ with target ranks $(1,\ell, \dots, \ell,1)$ assuming that $\ell = \Theta(R)$. We explain in some detail the analysis of the computational costs of the algorithms described in \Cref{ssec:ttrounding}, since it will help analyze the computational costs of the proposed algorithms.

\paragraph{Right-to-Left Orthogonalization (\Cref{alg:orthogonalizationr2l})}

The two main contributions to the leading order flop cost are: (1) in \cref{line:orthogonalizationr2l:qr}, the thin QR factorization of the matrix $\HOp(\TT{X}{n})$ of size $IR \times R$, which cost $4IR^{3} + \mathcal{O}(R^3)$ flops; (2) in \cref{line:orthogonalizationr2l:matmul}, the multiplication of a matrix of size $IR \times R$ with another matrix of size $R \times R$ that forms the product $\VOp(\TT{Y}{n-1})\M{R}^\top$. Since 
$\M{R}$ is triangular, the computation in \cref{line:orthogonalizationr2l:matmul} requires only $IR^{3}$ flops to leading order. 
Therefore, the resulting total cost of \Cref{alg:orthogonalizationr2l} is
\begin{eqnarray}\label{eq:cost_R2LO}
C_{\rm Orth}= I(N-2) \cdot 5R^3 + \mathcal{O}(IR^2 + NR^3)\quad \mbox{flops}.
\end{eqnarray}

\paragraph{TT-Rounding (\Cref{alg:TT-rounding-RLR})}

The main contributions to the total cost of TT-rounding come from two sources: in 
\cref{line:TTRound:OrthRL}, we first right-orthogonalize $\T{Y}$ to obtain $\T{X}{}$ using \Cref{alg:orthogonalizationr2l}; this cost has just been analyzed and is given by $C_{\rm Orth}$. Next we analyze the cost of the loop in 
\crefrange{line:TTRound:loop}{line:TTRound:Line10}.
In \cref{line:TTRound:loopQR}, we compute the QR factorization of matrix $\VOp(\TT{X}{n})$ of size $I\ell \times R$ using the Householder-QR algorithm, which costs $4I\ell R^{2} + \mathcal{O}(R^3)$ flops. Next, in \cref{line:TTRound:loopSVD} performing the SVD requires $\mathcal{O}(R^{3})$ flops.
The next two steps in lines 
\cref{line:TTRound:loopMatMul1,line:TTRound:loopMatMul2} prepare the factors for the next iteration.
The matrix-matrix multiplication in \cref{line:TTRound:loopMatMul1} costs $2IR\ell^2$.  
Notice that the size of $\TT{X}{n}$ has been reduced from the previous iteration, so that 
\cref{line:TTRound:loopMatMul2} requires $2IR^{2}\ell$ flops. Hence, the total cost of \Cref{alg:TT-rounding-RLR} is
\begin{align*}\label{eq: cost_TTR_ROLT}
C_{\rm TTR}=& \> C_{\rm Orth} + (N-2) \cdot \left(6 IR^{2}\ell + 2 IR\ell^2\right) + \mathcal{O}(IR^{2} + NR^3) \\ = & \> I(N-2) \cdot \left(5 R^{3} + 6R^{2}\ell + 2R\ell^2\right) + \mathcal{O}(IR^{2} + NR^3) \quad \mbox{flops}.
\end{align*}

\paragraph{Right-to-Left Partial Contraction (\Cref{alg:contractionr2l})}
Consider two TT-tensors $\T{X}$ and $\T{Y}$ with ranks $\{R_j^{\T{X}}\}$ and $\{R_j^{\T{Y}}\}$. For simplicity of exposition, we take $\ell = R_j^{\T{X}}$ and $R = R_j^{\T{Y}}$ for $j=1,\dots,N-1$; that is we take the same ranks in each core tensor.
We denote the computational cost of \Cref{alg:contractionr2l} 
by $C_{\rm Contr}$. The first contribution in the computational cost (\cref{line:contractionr2l:MatMul1,line:contractionr2l:MatMul2}) comes from two different matrix multiplications, i.e., matrix of size $I\ell \times \ell$ with the one of size $\ell \times R$ which costs  $2I \ell^2 R$ flops, and matrix of size $\ell \times IR$ with the one of size $I R \times R$ which costs $2IR^2 \ell$ flops. The cost of 
\cref{line:contractionr2l:MatMul0} is of order $\mathcal{O}(IR \ell)$ flops. 
The overall cost of \Cref{alg:contractionr2l} is, therefore,
\begin{eqnarray}\label{eq:cost_R2LC}
C_{\rm Contr}=I(N-2) \cdot ( 2R\ell^2 + 2R^2 \ell ) + \mathcal{O}(IR \ell) \quad \mbox{flops}.
\end{eqnarray}
{If $\ell =R$, then we are computing an inner product and the associated cost is $C_{\rm Contr}=4(N-2)IR^3 + \mathcal{O}(IR^2)$ flops.} 

\section{Additional Tensor network diagrams}\label{sec:networkdiagrams}

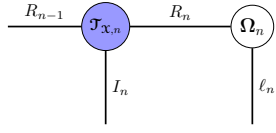
\begin{figure}[!ht]
    \centering
    \scalebox{0.65}{
     \tikzstyle{vertex} = [circle,draw,scale=.8]
\tikzstyle{edge} = [line width=1pt]

\begin{tikzpicture}[scale=1]

\node[vertex,fill=blue!40] (x) at (0,0) {\large $\TT{X}{n}$};

\node[vertex] (omega) at (3,0) {\Large $\M{\Omega}_n$};

\draw[edge] (x) -- (-2,0) node[midway,above] {$R_{n-1}$};
\draw[edge] (x) -- (omega) node[midway,above] {$R_{n}$};
\draw[edge] (x) -- (0,-2) node[midway,right] {$I_n$};
\draw[edge] (omega) -- (3,-2) node[midway,right] {$\ell_n$};

\end{tikzpicture}
    }
    \caption{Random projection steps for Orthogonalize-then-Randomize \Cref{alg:TT-rounding-Rand1}.}
    \label{fig:orth-rand}
\end{figure}

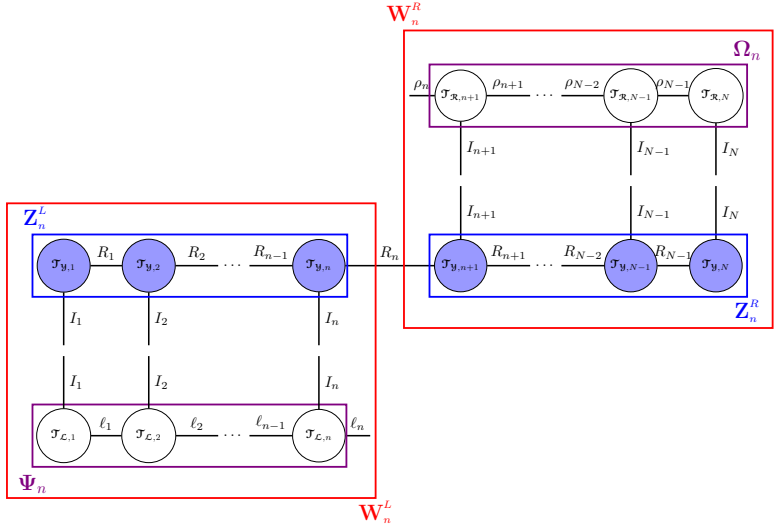
\begin{figure}[!ht]
    \centering
    \scalebox{0.58}{\tikzstyle{vertex} = [circle,draw,scale=.8]
\tikzstyle{edge} = [line width=1pt]

\begin{tikzpicture}[scale=.65]

\node[vertex,fill=blue!40] (a0) at (0,0) {$\hspace{0.095in}\TT{Y}{1\hspace{0.095in}}$};
\node[vertex,fill=blue!40] (a1) at (3,0) {$\hspace{0.095in}\TT{Y}{2}\hspace{0.095in}$};
\node (a2) at (6,0) {$\cdots$};
\node[vertex,fill=blue!40] (a3) at (9,0) {$\hspace{0.085in}\TT{Y}{n}\hspace{0.085in}$};
\node[vertex,fill=blue!40] (a4) at (14,0) {$\hspace{0.015in}\TT{Y}{n+1}\hspace{0.015in}$};
\node (a5) at (17,0) {$\cdots$};
\node[vertex,fill=blue!40] (a6) at (20,0) {$\TT{Y}{N-1}$};
\node[vertex,fill=blue!40] (a7) at (23,0) {$\hspace{0.07in}\TT{Y}{N}\hspace{0.07in}$};

\foreach \i in {0,...,3} {
    \node (b\i) at (\i*3,-3) {};
}
\foreach \i in {4,...,7} {
    \node (b\i) at (\i*3+2,3) {};
}


\node[vertex] (c0) at (0,-6) {$\hspace{0.095in}\TT{L}{1}\hspace{0.095in}$};
\node[vertex] (c1) at (3,-6) {$\hspace{0.095in}\TT{L}{2}\hspace{0.095in}$};
\node (c2) at (6,-6) {$\cdots$};
\node[vertex] (c3) at (9,-6) {$\hspace{0.08in}\TT{L}{n}\hspace{0.08in}$};
\node (c34a) at (11,-6) {};

\node (c34b) at (12,6) {};
\node[vertex] (c4) at (14,6) {$\TT{R}{n+1}$};
\node (c5) at (17,6) {$\cdots$};
\node[vertex] (c6) at (20,6) {$\TT{R}{N-1}$};
\node[vertex] (c7) at (23,6) {$\hspace{0.075in}\TT{R}{N}\hspace{0.075in}$};

\draw[edge] (a0) -- (a1) node[midway,above] {$R_1$};
\draw[edge] (a1) -- (a2) node[midway,above] {$R_2$};
\draw[edge] (a2) -- (a3) node[midway,above] {$R_{n-1}$};
\draw[edge] (a3) -- (a4) node[midway,above] {$R_{n}$};
\draw[edge] (a4) -- (a5) node[midway,above] {$R_{n+1}$};
\draw[edge] (a5) -- (a6) node[midway,above] {$R_{N-2}$};
\draw[edge] (a6) -- (a7) node[midway,above] {$R_{N-1}$};

\draw[edge] (c0) -- (c1) node[midway,above] {$\ell_1$};
\draw[edge] (c1) -- (c2) node[midway,above] {$\ell_2$};
\draw[edge] (c2) -- (c3) node[midway,above] {$\ell_{n-1}$};
\draw[edge] (c3) -- (c34a) node[midway,above] {$\ell_n$};
\draw[edge] (c34b) -- (c4) node[midway,above] {$\rho_n$};
\draw[edge] (c4) -- (c5) node[midway,above] {$\rho_{n+1}$};
\draw[edge] (c5) -- (c6) node[midway,above] {$\rho_{N-2}$};
\draw[edge] (c6) -- (c7) node[midway,above] {$\rho_{N-1}$};

\draw[edge] (a0) -- (b0) node[midway,right] {$I_1$};
\draw[edge] (a1) -- (b1) node[midway,right] {$I_2$};
\draw[edge] (a3) -- (b3) node[midway,right] {$I_{n}$};
\draw[edge] (a4) -- (b4) node[midway,right] {$I_{n+1}$};
\draw[edge] (a6) -- (b6) node[midway,right] {$I_{N-1}$};
\draw[edge] (a7) -- (b7) node[midway,right] {$I_{N}$};

\draw[red, very thick] (-2,2.2) rectangle (11,-8.2);
\node[text=red] at (11,-8.7) {$\Large \textbf{W}_n^L$};

\draw[blue, very thick] (-1.1,1.1) rectangle (10,-1.1);
\node[text=blue] at (-1,1.7) {$\Large \textbf{Z}_{n}^L$};

\draw[violet, very thick] (-1.1,-4.9) rectangle (9.97,-7.1);
\node[text=violet] at (-1.1,-7.7) {$\Large \M{\Psi}_n$};

 \draw[edge] (c0) -- (b0) node[midway,right] {$I_1$};
 \draw[edge] (c1) -- (b1) node[midway,right] {$I_2$};
 \draw[edge] (c3) -- (b3) node[midway,right] {$I_{n}$};
 \draw[edge] (c4) -- (b4) node[midway,right] {$I_{n+1}$};
 \draw[edge] (c6) -- (b6) node[midway,right] {$I_{N-1}$};
 \draw[edge] (c7) -- (b7) node[midway,right] {$I_{N}$};

\draw[red, very thick] (12,8.3) rectangle (25,-2.2);
\node[text=red] at (12,8.9) {$\Large \textbf{W}_n^R$};

\draw[violet, very thick] (12.9,7.1) rectangle (24.1,4.9);
\node[text=violet] at (24.1,7.6) {$\Large \M{\Omega}_n$};

\draw[blue, very thick] (12.9,1.1) rectangle (24.1,-1.1);
\node[text=blue] at (24.1,-1.6) {$\Large \textbf{Z}_{n}^R$};

\end{tikzpicture}}
    \caption{Random projections for Two-Sided-Randomization (Generalized Nystr\"{o}m) \Cref{alg:TT-rounding-Rand3}.}
    \label{fig:two-sided}
\end{figure}

\begin{figure}[!ht]
    \centering
    \scalebox{0.55}{\tikzstyle{vertex} = [circle,draw,scale=.8]
\tikzstyle{edge} = [line width=1pt]

\begin{tikzpicture}[scale=.6]

\node[vertex,fill=blue!40] (a0) at (0,0) {$\hspace{0.095in}\TT{Y}{1\hspace{0.095in}}$};
\node[vertex,fill=blue!40] (a1) at (3,0) {$\hspace{0.095in}\TT{Y}{2}\hspace{0.095in}$};
\node (a2) at (6,0) {$\cdots$};
\node[vertex,fill=blue!40] (a3) at (9,0) {$\hspace{0.085in}\TT{Y}{n}\hspace{0.085in}$};
\node[vertex,fill=blue!40] (a3a) at (12,0) {$\hspace{0.14in}\M{L}_{n}\hspace{0.14in}$};
\node[vertex,fill=blue!40] (a3b) at (15,0) {$\hspace{0.135in}\M{R}_{n}\hspace{0.135in}$};
\node[vertex,fill=blue!40] (a4) at (18,0) {$\hspace{0.015in}\TT{Y}{n+1}\hspace{0.015in}$};
\node (a5) at (21,0) {$\cdots$};
\node[vertex,fill=blue!40] (a6) at (24,0) {$\TT{Y}{N-1}$};
\node[vertex,fill=blue!40] (a7) at (27,0) {$\hspace{0.07in}\TT{Y}{N}\hspace{0.07in}$};

\foreach \i in {0,...,3} {
    \node (b\i) at (\i*3,-3) {};
}
\foreach \i in {4,...,7} {
    \node (b\i) at (\i*3+6,3) {};
}


\draw[edge] (a0) -- (a1) node[midway,above] {$R_1$};
\draw[edge] (a1) -- (a2) node[midway,above] {$R_2$};
\draw[edge] (a2) -- (a3) node[midway,above] {$R_{n-1}$};
\draw[edge] (a3) -- (a3a) node[midway,above] {$R_{n}$};
\draw[edge] (a3a) -- (a3b) node[midway,above] {$\ell_n$};
\draw[edge] (a3b) -- (a4) node[midway,above] {$R_{n}$};
\draw[edge] (a4) -- (a5) node[midway,above] {$R_{n+1}$};
\draw[edge] (a5) -- (a6) node[midway,above] {$R_{N-2}$};
\draw[edge] (a6) -- (a7) node[midway,above] {$R_{N-1}$};


\draw[edge] (a0) -- (b0) node[midway,right] {$I_1$};
\draw[edge] (a1) -- (b1) node[midway,right] {$I_2$};
\draw[edge] (a3) -- (b3) node[midway,right] {$I_{n}$};
\draw[edge] (a4) -- (b4) node[midway,right] {$I_{n+1}$};
\draw[edge] (a6) -- (b6) node[midway,right] {$I_{N-1}$};
\draw[edge] (a7) -- (b7) node[midway,right] {$I_{N}$};


\draw[blue, very thick] (-1.5,1.5) rectangle (10.5,-1.5);
\node[text=blue] at (-1,2.5) {$\Large \textbf{Z}_{n}^L$};





\draw[blue, very thick] (16.5,1.5) rectangle (28.5,-1.5);
 \node[text=blue] at (28,-2.5) {$\Large \textbf{Z}_{n}^R$};

\end{tikzpicture}}
    \caption{Internal left and right factors in the Two-Sided-Randomization (Generalized Nystr\"{o}m) \Cref{alg:TT-rounding-Rand3}.}
    \label{fig:two-sided1}
\end{figure}
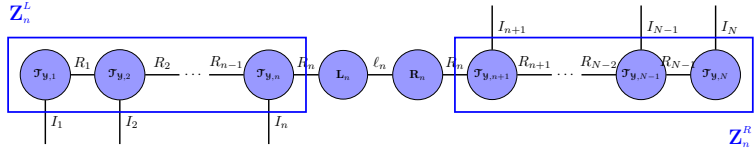

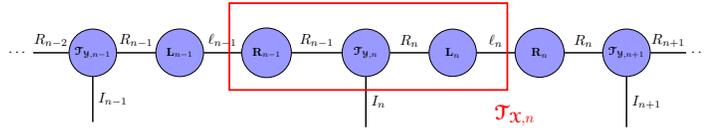
\begin{figure}[!ht]
    \centering
    \scalebox{0.55}{\tikzstyle{vertex} = [circle,draw,scale=.8]
\tikzstyle{edge} = [line width=1pt]

\begin{tikzpicture}[scale=.6]

\node (a1) at (0,0) {$\cdots$};
\node[vertex,fill=blue!40] (a2) at (3,0) {$\TT{Y}{n-1}$};
\node[vertex,fill=blue!40] (a3) at (6.5,0) {$\hspace{0.06in}\M{L}_{n-1}\hspace{0.06in}$};
\node[vertex,fill=blue!40] (a4) at (10,0) {$\hspace{0.06in}\M{R}_{n-1}\hspace{0.06in}$};
\node[vertex,fill=blue!40] (a5) at (14,0) {$\hspace{0.07in}\TT{Y}{n}\hspace{0.07in}$};
\node[vertex,fill=blue!40] (a6) at (17.5,0) {$\hspace{0.13in}\M{L}_{n}\hspace{0.13in}$};
\node[vertex,fill=blue!40] (a7) at (21,0) {$\hspace{0.13in}\M{R}_{n}\hspace{0.13in}$};
\node[vertex,fill=blue!40] (a8) at (24.5,0) {$\TT{Y}{n+1}$};
\node (a9) at (27.5,0) {$\cdots$};

\foreach \i in {0,...,3} {
     \node (b\i) at (\i*3,-3) {};
 }


 \draw[edge] (a1) -- (a2) node[midway,above] {$R_{n-2}$};
\draw[edge] (a2) -- (a3) node[midway,above] {$R_{n-1}$};
 \draw[edge] (a3) -- (a4) node[midway,above] {$\ell_{n-1}$};
 \draw[edge] (a4) -- (a5) node[midway,above] {$R_{n-1}$};
 \draw[edge] (a5) -- (a6) node[midway,above] {$R_{n}$};
 \draw[edge] (a6) -- (a7) node[midway,above] {$\ell_{n}$};
\draw[edge] (a7) -- (a8) node[midway,above] {$R_{n}$};
\draw[edge] (a8) -- (a9) node[midway,above] {$R_{n+1}$};

 \draw[edge] (a2) -- (b1) node[midway,right] {$I_{n-1}$};
 \draw[edge] (a5) -- (14,-3) node[midway,right] {$I_n$};
\draw[edge] (a8) -- (24.5,-3) node[midway,right] {$I_{n+1}$};

\draw[red, very thick] (8.5,2) rectangle (19.7,-1.5);
\node[text=red] at (20,-2.5) {$\Large \TT{X}{n}$};

\end{tikzpicture}}
    \caption{Left and right factors distributions within the Two-Sided-Randomization (Generalized Nystr\"{o}m) \Cref{alg:TT-rounding-Rand2}.}
    \label{fig:two-sided2}
\end{figure}

In \Cref{fig:orth-rand,fig:one-sided}, we give tensor network diagrams that help illustrate the random projection steps in \Cref{alg:TT-rounding-Rand1} and \Cref{alg:TT-rounding-Rand2}. In \Cref{fig:two-sided,fig:two-sided1,fig:two-sided2}, we provide the tensor network diagrams illustrating \Cref{alg:TT-rounding-Rand3}.

\section{Additional Numerical Experiments: Hilbert-type TT-tensor}\label{sec:hilbert}

To supplement the numerical experiments in the main manuscript, we propose a Hilbert-type TT-tensor, i.e., a TT-tensor $\T{X} \in \mathbb{R}^{I_1 \times \dots \times I_N}$ with the cores of the form
\[ 
\TT{X}{n} = \frac{1}{i_1 + i_2 + i_3 -1 }, \qquad 1 \leq i_1 \leq r_{n-1}, \> 1 \leq i_2 \leq I_n, \> 1 \leq i_3 \leq r_n, 
\]
for $n=1,\dots,N$. Our definition of a Hilbert-type TT-tensor is inspired by Hilbert matrices, i.e., symmetric and positive definite matrices with rapidly decaying singular values. Since the singular values of each mode-unfolding of the cores of the Hilbert-type TT-tensor decay very quickly, it is a useful tensor for validating the accuracy of the randomized algorithms proposed in this paper. 

\begin{figure}[h!]
    \centering
    \includegraphics[scale=0.6]{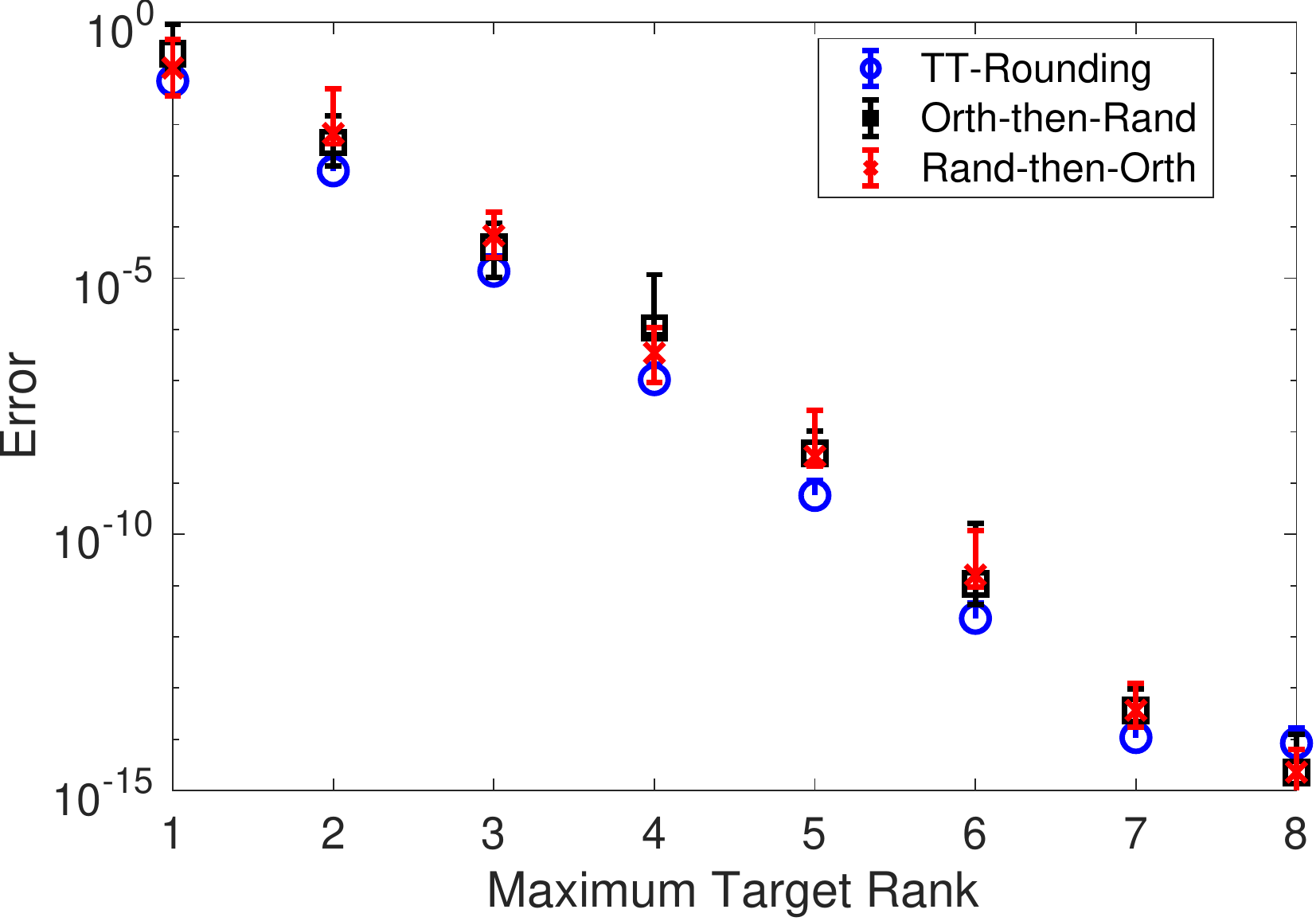}
    \caption{Hilbert TT-tensor truncation. Relative Error between a low rank approximation obtained via deterministic or randomized TT-rounding, and a full rank Hilbert TT-tensor.}
    \label{Fig:Varyingp}
\end{figure}

In our test, we consider the order-$10$ Hilbert-type TT-tensor with dimensions $I_1 = 6, I_2 = 8, \ldots, I_{10} = 24$ and ranks $(1,4,5,6,\ldots,11,12,1)$ and compare numerically the deterministic approach (\Cref{alg:TT-rounding-RLR}), the Orthogonalize-then-Randomize algorithm and the  Randomize-then-Orthogonalize algorithm (\Cref{alg:TT-rounding-Rand1,alg:TT-rounding-Rand2}, respectively). Results of our experiments are presented in \Cref{Fig:Varyingp}, with all three algorithms exhibiting a very similar accuracy.

\section{Additional details regarding \Cref{sec:PDE_Example_experiments}}
\label{sec:PDE_Example}

\begin{figure}[hb!]
    \centering
    \includegraphics[width=.5\textwidth]{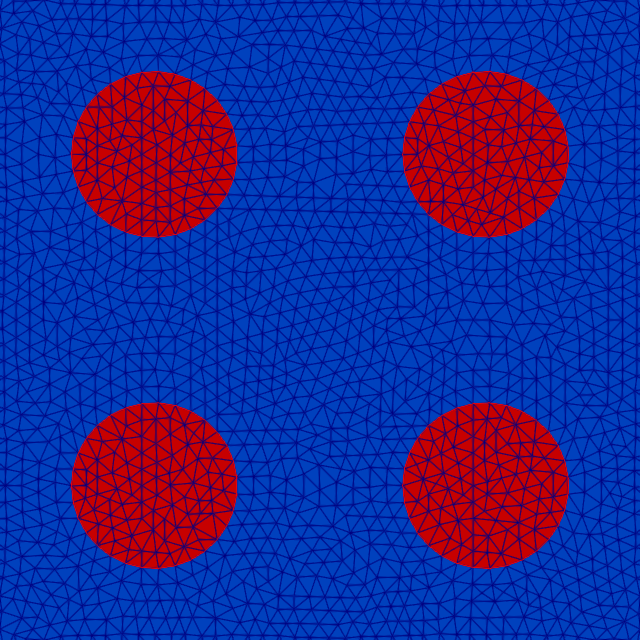}
    \caption{Geometry and mesh used for the cookie problem~\eqref{eq:cookie} with $N=5$.}
    \label{Fig:cookies}
\end{figure}

In this section, we present additional details regarding the \emph{cookie problem} we examined in \Cref{sec:PDE_Example_experiments}.

\begin{equation}
\begin{aligned}
  \label{eq:cookie}
  -\text{div}(\sigma(x,y;{\bm \rho}) \nabla(u(x,y; {\bm \rho}))) &= f(x,y) \quad &\text{ in } \Omega,\\
   u(x,y; {\bm \rho}) &= 0 \quad &\text{ on } \partial \Omega,
\end{aligned}
\end{equation}
where $\Omega$ is $(-1,1)\times (-1,1)$, $\partial \Omega$ is the boundary of $\Omega$ and $\sigma$ is defined as follows
\[
  \sigma(x,y;{\bm \rho}) =
        \begin{cases}
                 1 + \rho_i & \text{if } (x,y) \in D_i,\\ 
                 1          & \text{elsewhere},
        \end{cases}
\]
where $D_i$ for $i=1,\ldots,N-1$ are disjoint disks distributed in $\Omega$ such that their centers are equidistant (see \Cref{Fig:cookies}) and $\rho_i$ is selected from a discrete set of $I_{i+1}$ samples $J_i \subset \R$ for $i = 1, \ldots,N-1$.

\noindent
Discretizing this problem in $\Omega$ yields the set of $\prod_{i=1}^{N-1} I_{i+1}$ parametrized linear systems of equations
\[
\left(\M{A}_{1,1} + \sum_{i=1}^{N-1}\rho_i \M{A}_{i+1,1}\right) {\bm x} = {\bm f},
\]
where $\M{A}_{1,1}\in\R^{I_1 \times I_1}$ is the discretization of the operator 
$-\text{div}(\sigma(x,y;{\bm \rho}) \nabla(u(x,y; {\bm \rho})))$ in $\Omega$ and $\M{A}_{i+1,1} \in\R^{I_1 \times I_1}$ are the discretizations of $-\text{div}(\chi_{D_i}\nabla(\cdot))$ in $\Omega$ where $\chi_{S}$ denotes the indicator function of the set $S$, and ${\bm f}$ is the discretization of the function $f$.

\bibliographystyle{abbrv}
\bibliography{refs}

\end{document}